\documentclass{amsart}
\usepackage[english]{babel}
\usepackage[T1]{fontenc}
\usepackage{geometry}
\usepackage[style=alphabetic,natbib=true,giveninits=true,url=false,maxbibnames=99]{biblatex}
\usepackage{hyperref}

\usepackage{graphicx} 
\usepackage{float}
\usepackage{braket}
\usepackage{amsmath}
\usepackage{amsfonts}
\usepackage{amssymb}
\usepackage{mathtools}
\usepackage[onehalfspacing]{setspace}
\usepackage{stmaryrd}
\usepackage{enumitem}

\usepackage{xcolor}

\newtheorem{theorem}{Theorem}[section]
\newtheorem{corollary}[theorem]{Corollary}
\newtheorem{lemma}[theorem]{Lemma}
\newtheorem{conjecture}[theorem]{Conjecture}
\newtheorem{prop}[theorem]{Proposition}

\theoremstyle{definition}
\newtheorem{definition}[theorem]{Definition}
\newtheorem{remark}[theorem]{Remark}

\newcommand{\Z}{\mathbb{Z}}
\newcommand{\N}{\mathbb{N}}
\newcommand{\Q}{\mathbb{Q}}
\newcommand{\R}{\mathbb{R}}

\newcommand{\fa}{\forall}

\DeclareMathOperator{\dimH}{\mathrm{dim}_{\mathrm{H}}}

\title{New portions of $M\setminus L$ and a lower bound on the Hausdorff distance between $L$ and $M$}

\begin{document}

\author{Clément Rieutord}
\address[Clément Rieutord]{Ecole Polytechnique, Rte de Saclay, 91120 Palaiseau}
\email{clement.rieutord@polytechnique.edu}

\author{Carlos Gustavo Moreira}
\address[Carlos Gustavo Moreira]{SUSTech International Center for Mathematics, Shenzhen, Guangdong, People’s Republic of China; \hfill\break
IMPA, Estrada Dona Castorina 110, 22460-320, Rio de Janeiro, Brazil}
\email{gugu@impa.br}

\author{Harold Erazo}
\address[Harold Erazo]{IMPA, Estrada Dona Castorina 110, 22460-320, Rio de Janeiro, Brazil}
\email{harold.eraz@gmail.com}

\thanks{The first author is partially supported by the Ecole Polytechnique, the second author is partially supported by CNPq and FAPERJ. The third author is partially supported by CAPES and FAPERJ}


\begin{abstract}
	Let $M$ and $L$ be the Markov and Lagrange spectra, respectively. It is known that $L$ is contained in $M$ and Freiman showed in 1968 that $M\setminus L\neq \emptyset$. In 2018 the first region of $M\setminus L$ above $\sqrt{12}$ was discovered by C. Matheus and C. G. Moreira, thus disproving a conjecture of Cusick of 1975. In 2022, the same authors together with L. Jeffreys discovered a new region near 3.938. In this paper, we will study two new regions of $M\setminus L$ above $\sqrt{12}$, in the vicinity of the Markov value of two periodic words of odd length that are non semisymmetric, which are $\overline{212332111}$ and $\overline{123332112}$. We will demonstrate that for both cases, there is a maximal gap of $L$ and a Gauss-Cantor set inside this gap that is contained in $M$. Moreover we show that at the right endpoint of those gaps we have local Hausdorff dimension equal to $1$.
	
	After studying the mentioned examples, we will provide a lower bound for the value of $d_H(M,L)$ (the Hausdorff distance between $M$ and $L$). 
\end{abstract}

\maketitle
\tableofcontents

\section{Introduction}

\begin{definition}
The best constant of Diophantine approximation $k(\alpha)$ for $\alpha \in \R \setminus \Q$ is given by:
$$l(\alpha) =\limsup
\limits_{\substack{p,q \rightarrow \infty\\
p,q \in \Z^2}}\frac{1}{|q(q\alpha -p)|}.$$
    The Lagrange spectrum is:
    $$L = \left\{\,l(\alpha) \:\:/ \:\: \alpha \in \R \setminus \Q\:\:\text{and}\:\: l(\alpha)< \infty\right\}.$$

\end{definition}

\begin{definition}
    Let $q(x,y) = ax^2+bxy+cy^2$ be a binary quadratic form with $(a,b,c)\in \R^3$ and $\Delta(q) = b^2-4ac$. We say $q$ is indefinite if $\Delta(q)>0$. 
    
    The Markov Spectrum is:
    $$M = \left\{\frac{\sqrt{\Delta(q)}}{\inf_{(x,y)\in \Z^2\setminus\{(0,0)\}}q(x,y)}\:\:/\:\: q \:\:\text{is an indefinite binary quadratic form}\right\}.$$

\end{definition}

The Lagrange and Markov spectra were first studied by Markov around 1880 (see \cite{Markov}). In particular he established that $L\cap(0,3)=M\cap(0,3)$ is a discrete set that accumulates at 3. M. Hall showed \cite{Hallray} in 1947 that $[c,\infty)\subset L$ for some positive constant $c>3$. In \cite{Freiman73} and \cite{Schecker77}, it was proved that in fact $[\sqrt{21},\infty)\subset L$ and finally in 1975, Freiman \cite{Freimanconstant} determined the largest half-ray contained in $L$, namely $[\mu,\infty)\subset L$, where $\mu = 4.527829\dots$ is an explicit quadratic surd. The set $[\mu,\infty)$ is known as \emph{Hall's ray}.

It is well known that $L\subset M$ are closed subsets of $\R$ (see \cite{Cusick-Flahive}). It took some time to decide if they were equal: only in 1968 Freiman \cite{F:non-coincidence} exhibited a countable subset of isolated points in $M\setminus L$ near 3.11; after that, Freiman proved in 1973 that $M\setminus L$ contains a point $\alpha_\infty$ near 3.29, and Flahive \cite{Fl77} showed in 1977 that $\alpha_\infty$ is the accumulation point of a countable subset of $M\setminus L$. These results led Cusick \cite{Cusick75} to conjecture in 1975 that the Lagrange and Markov spectra coincide after $\sqrt{12}$, i.e. $(M\setminus L)\cap[\sqrt{12},\infty)=\varnothing$.

Only recently our knowledge of $M\setminus L$ changed significantly: in \cite{FreimanExample}, Carlos Matheus and Carlos Gustavo Moreira showed that in fact there is a Gauss-Cantor set inside $M\setminus L$ that contains $\alpha_\infty$ (in particular $M\setminus L$ is uncountable).\footnote{In the PhD thesis by \cite{Michiganthesis} in 1976, the author found that there is a Gauss-Cantor set near $\alpha_\infty$ contained in $M\setminus L$. Consequently, he proved that $M\setminus L$ is uncountable. However, it was not observed that this already implies that $M\setminus L$ has positive Hausdorff dimension, nor was there an attempt made to completely characterize the region. It is unfortunate that this work went essentially unnoticed; it is not even mentioned in the book of \cite{Cusick-Flahive}.} As a by product of their work they showed that $\dimH(M\setminus L)>0.353$, where $\dimH(M\setminus L)$ denotes the Hausdorff dimension of the set $M\setminus L$. Shortly after they showed that Cusick conjecture mentioned above is false: in \cite{MM:markov_lagrange} they found a Gauss-Cantor set near 3.7 of Hausdorff dimension larger than 1/2. Moreover they showed that in fact $\dimH(M\setminus L)<0.987$ (in particular $L$ and $M$ have the same interior). The same authors together with  M. Pollicott and P. Vytnova proved rigorously in \cite{GMPP} that
\begin{equation*}
    0.537152 < \dimH(M \setminus L) < 0.796445.
\end{equation*}
The upper bound is the current best upper bound on the Hausdorff dimension of $M\backslash L$, while in \cite{LukeGuguMatheus} the lower bound was improved to 0.593.

By a \emph{region} of $M\setminus L$ we mean a non-empty intersection of $M$ with a maximal gap of $L$. In \cite{Freiman311} the region near 3.11 related to Freiman example of 1968 was completely characterized: again the authors found a Gauss-Cantor set contained on $M\setminus L$ (see the book \cite{CDMLS} for more details). 

In all the above examples, the regions of $M\setminus L$ turned out to be closed subsets of $\R$, more precisely, they consists of a Gauss-Cantor set together with an infinite set of isolated points. Motivated by this, T. Bousch asked whether the set $M\setminus L$ is a closed subset of $\R$. In a first attempt to answer negatively this question, in \cite{MminusLnear3} Lima-Matheus-Moreira-Vieira gave some evidence towards the possibility that $3\in L\cap\overline{(M\setminus L)}$, more precisely they exhibited a decreasing sequence of elements on $M$ and proved that the first four of them belong to $M\setminus L$. In spite of the fact that $3\in L\cap\overline{(M\setminus L)}$ remains an open question, the same authors managed to prove in \cite{MminusLisnotclosed} that $1+\frac{3}{\sqrt{2}}\in L\cap\overline{(M\setminus L)}$, by constructing a decreasing sequence of elements of $M\backslash L$ that converge to $1+\frac{3}{\sqrt{2}}\in L$.

More recently, in \cite{LukeGuguMatheus}, one new region of $M\setminus L$ above $\sqrt{12}$ was discovered near 3.938. This region gives the current known maximal elements of $M\setminus L$. In fact, by means of computational search, four new regions of $M\setminus L$ above $\sqrt{12}$ were suggested in \cite{LukeGuguMatheus}. In this paper we study one of the regions they suggested, namely the one associated with 123332112. The other region we study is the one associated with 212332111, which was not noticed before.

We say that a finite word $w$ is \emph{semisymmetric}\footnote{The term was coined by \cite{Fl77}, however this concept already appeared in the work of Berstein in \cite[Page 42]{Berstein73}} if it is a palindrome or the concatenation of two palindromes, equivalently, if the orbits by the shift map of the bi-infinite sequences $\overline{w}$ and $\overline{w^T}$ are equal. For example 212332111 and 123332112 are non semisymmetric words of odd length. As is explained in \cite{MM:markov_lagrange}, this is one of the key notions behind the construction of all known regions of $M\setminus L$.

For more details about the spectra, we refer the reader to the survey \cite{Spectrasurvey}, the classical textbook of Cusick-Flahive \cite{Cusick-Flahive}
or the book of Lima-Matheus-Moreira-Roma\~{n}a \cite{CDMLS} for a more modern point of view.

\subsection{Perron's description of the Lagrange and Markov spectra}
Perron managed to characterize the Lagrange and Markov spectra in terms of Dynamical Systems.

\begin{definition}
    Let us call $\Sigma = (\N^*)^{\Z}$ and $\lambda_0 : \Sigma \rightarrow \R$ the function:
    $$\fa \: (a_n)_{n\in \Z}, \:\: \lambda_0((a_n)_{n\in \Z}) = \left[a_0,a_1,a_2,\dots\right] + \left[0,a_{-1},a_{-2},\dots\right],$$
    with 
    $\left[a_0,a_1,a_2,\dots\right] = a_0+\frac{1}{a_1+\frac{1}{a_2+\dots}}$ the value of the continued fraction associated with the sequence $(a_n)_{n\geq 0}$.
    We also define the function $\lambda_k$, for $k\in \Z$ as follows: $\fa (a_n)_{n\in \Z}$, $\lambda_k((a_n)_{n\in \Z}) = \lambda_0((a_{n+k})_{n\in \Z})$.
\end{definition}
\begin{theorem}[Perron, 1921]
We define the Langrage and Markov value of a sequence $S \in (\N^*)^{\Z}$ by: 
$$l(S) = \limsup_{k\rightarrow \infty}\lambda_k(S) \:\:\text{and}\:\:  m(S) = \sup_{k \in \Z}\lambda_k(S).$$
We have:
$$L = \left\{l(S) \:\:/\:\: S \in (\N^*)^{\Z} \:\:\text{and}\:\: l(S)<\infty \right\} \:\:\text{and}\:\: M = \left\{m(S) \:\:/\:\: S \in (\N^*)^{\Z} \:\:\text{and}\:\: m(S)<\infty \right\}.$$

\end{theorem}
For the future, a sequence marked with a "$^*$" is a sequence where the index $0$ is next to the left of the asterisk.
For example, we have:
$$\lambda_0(\dots1223\textcolor{cyan}{3}^*211112\dots) = \left[\textcolor{cyan}{3},2,1,1,1,1,2,\dots\right] + \left[0,3,2,2,1,\dots\right].$$

\subsection{New portions of $M\setminus L$}

Our main result is the complete characterization of two new regions of $M\setminus L$. Moreover, we also give information about the local dimension at the right endpoint of the gaps of $L$ where these regions are contained.

\begin{theorem} \label{principal1}

    \item There are two non-empty open intervals $A_1$ and $ A_2$  such that:
    $$\fa i\in \{1,2\}, \:\: A_i\cap L = \emptyset \:\:\text{and} \:\:\exists \,K_i \:\:\text{a Gauss-Cantor set such that} \:\: K_i\subset A_i \cap (M \setminus L).$$
    Furthermore, we have: 
    \begin{enumerate}
        \item \begin{enumerate}
            \item $\inf A_1 = \lambda_0(\overline{21233^*2111}) \approx 3.6767$,
            \item $|A_1| \approx 8.42651\times 10^{-12}$,
            \item $\sup A_1  \in L'$,
            \item $\fa \alpha >0,\: \dimH\{M\cap (\sup A_1, \sup A_1+\alpha)\} = \dimH\{L\cap (\sup A_1, \sup A_1+\alpha)\} = 1$,
        \end{enumerate}
        
        \item \begin{enumerate}
            \item $\inf A_2 = \lambda_0(\overline{12333^*2112}) \approx 3.72627$,
            \item $|A_2| \approx 5.88429\times 10^{-12}$,
            \item $\sup A_2 \in L'$,
            \item $\fa \alpha >0,\: \dimH\{M\cap (\sup A_2, \sup A_2+\alpha)\} = \dimH\{L\cap (\sup A_2, \sup A_2+\alpha)\} = 1$,
        \end{enumerate}  
        where $L^\prime$ denotes the accumulation points of $L$.
    \end{enumerate}
    
\end{theorem}

For a detailed description of the regions of Theorem \ref{principal1} see Theorem \ref{description globale M/L} and Theorem \ref{description globale M/L second region}.

\subsection{Hausdorff distance between $M$ and $L$}

Another interesting way of measure how different the spectra are is by computing the \emph{Hausdorff distance} between them. Since $L\subset M$, we have that $d_H(L,M)=\sup_{m\in M\setminus L}d(L,m)$.

\begin{theorem}\label{principal2}
    The Hausdorff distance between $M$ and $L$ is at least $\delta_0$, where
    \begin{align*}
        \delta_0 &= \frac{272052036746460995 - 3973474319367040 \sqrt{87} - 2762049221999040 \sqrt{18229} + 353887557067187 \sqrt{151905}}{226488036203921280} \\ 
        &\approx 9.1094243388\times 10^{-8}.
    \end{align*}
\end{theorem}

The reason why the bound is significant, is because in all regions discovered so far, including the ones described in Theorem \ref{principal1}, the one that is contained in the largest gap of $L$ is the one given by the odd non semisymmetric word $2112122$ of length 7. This is precisely the region that contains the example discovered by Freiman in 1973 and that was later characterized in \cite{FreimanExample}. The intuition is that longer words should give regions of $M\setminus L$ that lie in smaller gaps. There are no odd non semisymmetric words in the alphabet $\{1,2\}$ of length less than 7. The shortest odd non semisymmetric word that contains a 3 and which is known to give a region of $M\setminus L$, was the region discovered in \cite{MM:markov_lagrange} which is associated to 1233222. The gap of $L$ that contains this region is roughly of the size $7\times 10^{-11}$ (the right border was corrected in \cite{Corrigendum}).

There are no known examples of $M\setminus L$ associated with odd non semisymmetric words that contain at least a 4. In fact there are no known examples above 3.945. For odd non semisymmetric words in the alphabet $\{1,2,3,4\}$ with at least one 3 or 4 with length less than 7, the method does not seem to work, because either the Markov value of the periodic word lies in Hall's ray or the property of local uniqueness is not know to hold (it could happen that the Markov value of these periodic words belongs to the interior of the spectra, so the local uniqueness could be impossible to establish). 

\begin{conjecture}
The Hausdorff distance between $M$ and $L$ is precisely $\delta_0$.   
\end{conjecture}

Most of this paper will consist in proving Theorem \ref{principal1}. In Section 2, we describe precisely the structure of $M \setminus L$ nearby $\lambda_0(\overline{21233^*211})$.
In Section 3, we describe  the structure of $M \setminus L$ nearby $\lambda_0(\overline{12333^*2112})$.
Finally, in Section 4, we give a proof of Theorem \ref{principal2}.

\section{Preliminaries}

In this paper, we will only work in the region $M\cap (-\infty,4)$. Since $m(\dots4\dots)> 4$, we see that we must only deal with sequences $S \in \{1,2,3\}^{\Z}$.

\begin{definition}
For a finite word $w\in (\N^*)^{(\Z)}$, we denote:
$w^n =w\dots w \:$ ($n$ times).  

We denote the bi-infinite periodic sequence: $\dots w\dots w\dots$ $=\overline{w}$.

Eventually, if $S$ is a bi-infinite sequence with a periodic side $w$ on the left (or on the right), we also write $S = \overline{w}\,S_{n}S_{n+1}\dots$ 
\end{definition}
\begin{remark}
We have for every finite word $w$:
$$m(\overline{w}) = l(\overline{w}).$$
And we write, for $w = w_1\dots w_n$:
$$\lambda_0(\overline{w_1\dots w_k^*\dots w_n}) = \lambda_0(\overline{w}\,w_1\dots w_k^*\dots w_n\,\overline{w}).$$
\end{remark}

\begin{definition}
    For all sequences $S = (S_n)_{n\in\Z}$, we define the transpose of $S$ by $S^T = (S_{-n})_{n\in\Z}$.
    
    Clearly we have $\lambda_0(S) = \lambda_0(S^T)$ and therefore $m(S) = m(S^T)$ and $l(S) = l(S^T)$.
\end{definition}

\begin{remark}
   In order to simplify notation, if $w$ is a finite word and $S$ a bi-infinite sequence, we write $S = w_{-p}\dots w_0^*\dots w_q$ to denote:
   $$S_{-p}\dots S_{0}\dots S_{q} = w.$$
   Here, the asterisk on the finite word represents the position $0$ on the sequence $S$.
\end{remark}

\begin{definition} 
Given $B=\{\beta_1,\dots,\beta_\ell\}$, $\ell\geq 2$, a finite alphabet of finite words $\beta_j\in(\mathbb{N}^*)^{r_j}$, $r_j\geq 1$ which is primitive (in the sense that $\beta_i$ does not begin by $\beta_j$ for all $i\neq j$) and a finite word $c_1\dots c_r\in(\mathbb{N}^*)^{r}$, $r\geq 0$, then the set $K\subseteq [0,1]$ defined by
\[K=\{[0;c_1,\dots,c_r,\gamma_1, \gamma_2, \dots] \ \mid\ \gamma_i\in B\},\]
is a \emph{Gauss-Cantor set} associated with $B$.
\end{definition}

For the sequel, we need the following property on the function $\lambda_0$ (see \cite[Lemma A.1]{geometricproperties} and \cite[Chapter 1, Lemma 2]{Cusick-Flahive}): 
\begin{prop}\label{continuite lambda}   
\begin{enumerate}(Properties of continuity of $\lambda_0$)
    \item $\lambda_0(R^T\,.)$ is bi-Lipschitz:
    
    $\fa\, R = R_1R_2\dots  \in (\N^*)^{\N^*}$,
    $ \fa S \in (\N^*)^{\N} $, we define the bi-infinite sequence $R^T\,S$ by $(R^T\,S)_{n} = S_n$ if $n \geq 0$ and $(R^T\,S)_{n} = R_{-n}$ if $n <0$. There are $C_1,C_2 >0$ such that $\fa S,S' \in \{1,2,3\}^{\N},\fa\, R = R_1R_2\dots  \in (\N^*)^{\N^*}$ we have:
    $$\frac{C_1}{9^N}\leq \left|\lambda_0(R^T\,S) - \lambda_0(R^T\,S')\right|\leq \frac{C_2}{2^N}$$
    with $N = \max\left\{n \:/ \: \fa \:0 \leq k\leq n,\: S_k = S_k'\right\}$.
    \item  $\lambda_0$ is Lipschitz:
    
    There is $C >0$ such that $\fa S,S' \in \{1,2,3\}^{\Z}$ we have:
    $$\left|\lambda_0(S) - \lambda_0(S')\right|\leq \frac{C}{2^N}$$
    with $N = \max\left\{n \:/ \: \fa |k|\leq n,\: S_k = S_k'\right\}$.
\end{enumerate}

\end{prop}

\subsection{Boundaries of finite words}
In order to estimate  the $\lambda_0$ value of finite words that only contain $\{1,2,3\}$, we will use the following bounding rules:
\begin{lemma}\label{encadrement13}
\begin{enumerate}
    \item $\lambda_0(\overline{13}\,a_{2n}\dots a_0^*\dots a_{2m}\overline{31}\,) \leq \lambda_0(\dots a_{2n}\dots a_0^*\dots a_{2m}\dots ) \leq\lambda_0(\overline{31}\,a_{2n}\dots a_0^*\dots a_{2m}\overline{13}\,)$,
    \item $\lambda_0(\overline{13}\,a_{2n}\dots a_0^*\dots a_{2m+1}\overline{13}\,) \leq \lambda_0(\dots a_{2n}\dots a_0^*\dots a_{2m+1}\dots ) \leq\lambda_0(\overline{31}\,a_{2n}\dots a_0^*\dots a_{2m+1}\overline{31}\,)$,
    \item $\lambda_0(\overline{31}\,a_{2n+1}\dots a_0^*\dots a_{2m}\overline{31}\,) \leq \lambda_0(\dots a_{2n+1}\dots a_0^*\dots a_{2m}\dots ) \leq\lambda_0(\overline{13}\,a_{2n+1}\dots a_0^*\dots a_{2m}\overline{13}\,)$,
    \item $\lambda_0(\overline{31}\,a_{2n+1}\dots a_0^*\dots a_{2m+1}\overline{13}\,) \leq \lambda_0(\dots a_{2n+1}\dots a_0^*\dots a_{2m+1}\dots ) \leq\lambda_0(\overline{13}\,a_{2n+1}\dots a_0^*\dots a_{2m+1}\overline{31}\,)$.
\end{enumerate}
\end{lemma}

\begin{figure}[H]
\centerline{\includegraphics[scale = 0.6]{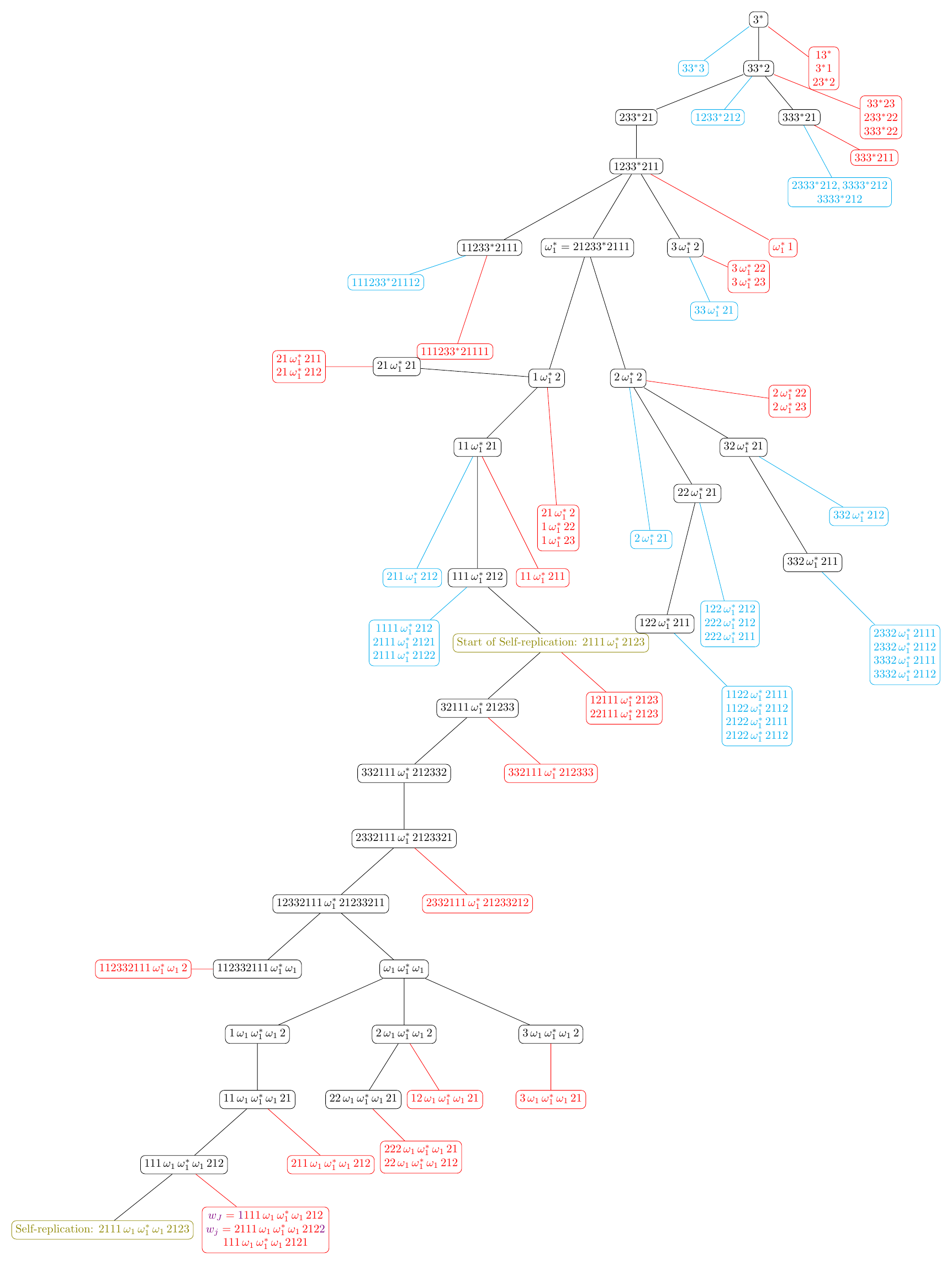}}
\caption{Sequence development Tree of $\overline{21233^*2111}$}
\label{Tree associated to omega_1}
\end{figure}

\begin{figure}[H]
\centerline{\includegraphics[scale = 0.6]{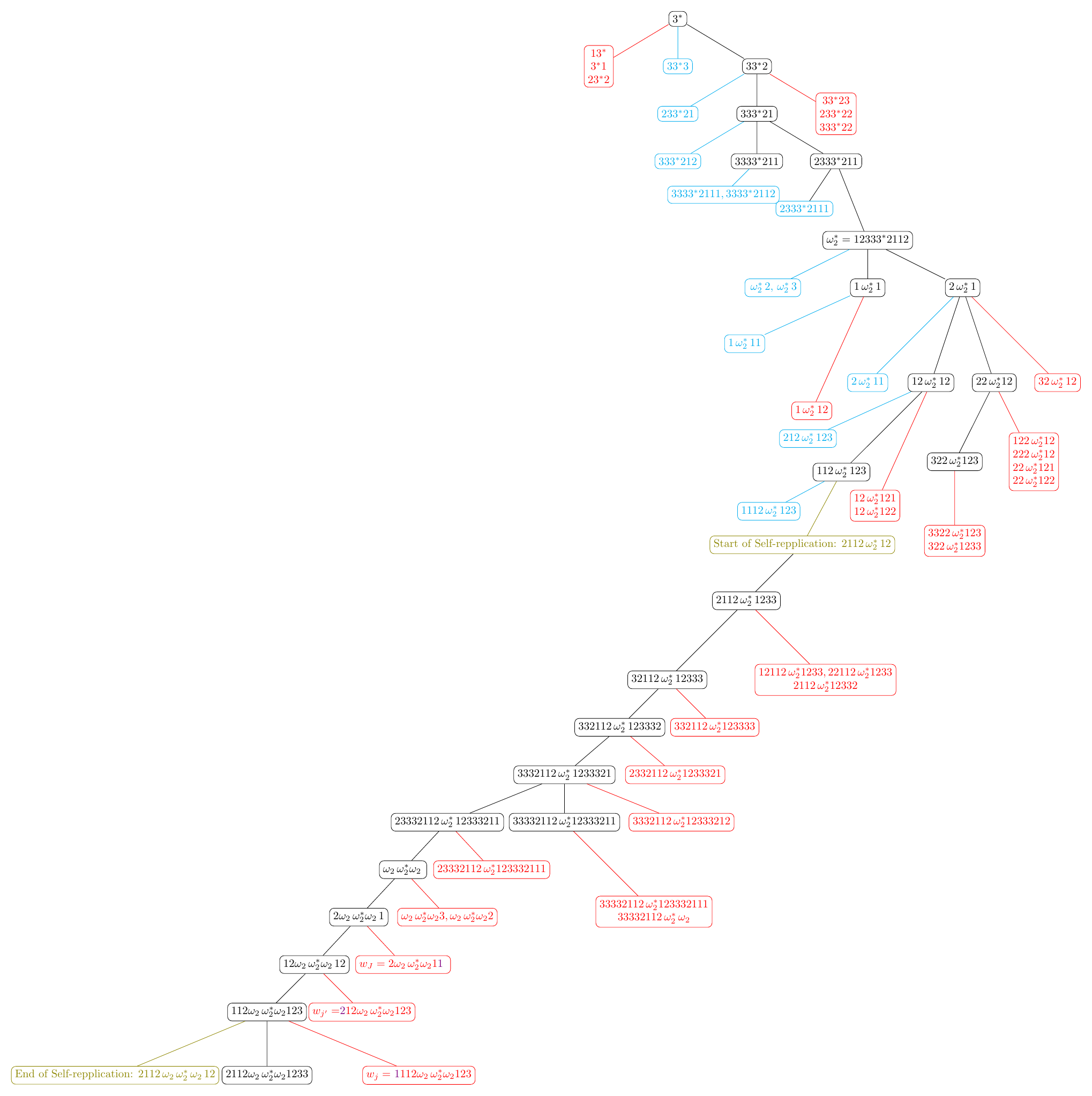}}
\caption{Sequence development Tree of $\overline{12333^*2112}$}
\label{Tree associated to w}
\end{figure}

These boundaries are a little gross, because the word "13" is quickly going to be forbidden, but it is enough to construct an accurate tree of possibilities (Figure \ref{Tree associated to omega_1} and \ref{Tree associated to w}) and obtain the results we want.

\subsection{Forbidden words}
In the following pages, we refer a lot to forbidden words. In each section we fix a real number $j_0$.
\begin{definition}
    A forbidden word is a finite word $w =a_{-n}\dots a_0\dots a_{m} \in \{1,2,3\}^{(\N)}$ that verifies 
$$\inf\{m(S)\:/ \:S\in\{1,2,3\}^\Z, w\text{ is subword of } S\} >j_0+\eta, \: \text{for some} \: \eta>0.$$
In other terms, $w$ is forbidden if for any bi-infinite sequence $S\in\{1,2,3\}^\Z$ containing $w$ as subword, the Markov value of $S$ is at least $j_0+\eta$ for some $\eta>0$.

\end{definition}

In practice we will forbid words according to the bounding rules of Lemma \ref{encadrement13}:

$$\lambda_0(\dots w\dots )\geq\min_{\varepsilon_l,\varepsilon_r \in \{\overline{13},\overline{31}\}}\lambda_0(\varepsilon_l\,w\,\varepsilon_r)  $$

\begin{remark}
    Because $\lambda_0(S) = \lambda_0(S^T)$, if $w$ is forbidden, then we also have that $w^T$ is forbidden as well.
\end{remark}

\section{Portion of $M\setminus L$ in the vicinity of $\lambda_0(\overline{21233^*2111})$}
In this section, we study the structure of $M$ in the vicinity of $j_0 =  \lambda_0(\overline{21233^*2111}) \approx 3.6766994172$.

Our aim is to determine the largest value $j_1$ for which $L\cap (j_0,j_1)= \emptyset$ is true, investigate the fractal structure of $(M\setminus L) \cap (j_0,j_1)$ and provide a description of the local structure of $L$ above $j_1$. To accomplish this, we will establish two key properties of the function $\lambda_0$ on this region:
\begin{enumerate}
    \item Local uniqueness. 
    \item Self-replication.
\end{enumerate}
We will resume the process of local uniqueness and self-replication in a tree (Figure \ref{Tree associated to omega_1}).

 We call for the future $\omega_1$ the finite word $212332111$. 
 We also denote $\omega_1 ^*$ the word $21233^*2111$, where the asterisk represents the position zero.

\subsection{Local uniqueness}
In this part, we will show that if $S$ verifies $|\lambda_0(S)-j_0|<\varepsilon$, with $ \varepsilon \approx 10^{-6}, $ then $S$ must be writen in the following form:
$$S_{-8}\dots S_{8} = 2111\; \omega_1^*\;  2123.$$
We will also determinate all forbidden words useful for the self-replication.

\begin{lemma}\label{fw1_1}(Forbidden words I):\\
Let $S \in \{1,2,3\}^{\Z}$.
\begin{enumerate}
    \item If $S = 13^*$, then $\lambda_0(S)> j_0+10^{-1}$.
    \item If $S = 23^*2$, then $\lambda_0(S)> j_0+10^{-2}$.
    \item If $S = 33^*23$, then $\lambda_0(S)> j_0+10^{-2}$.
    \item If $S = 233^*22 \:\:\text{or} \:\:333^*22 $, then $\lambda_0(S)> j_0+10^{-2}$.
    \item If $S = 333^*211$, then $\lambda_0(S)> j_0+10^{-3}$.
    \item If $S = 1233^*2112$, then $\lambda_0(S)> j_0+10^{-3}$.
    \item If $S = 111233^*21111$, then $\lambda_0(S)> j_0+10^{-4}$.
\end{enumerate}
\end{lemma}

As a consequence, if $S$ verifies $m(S)\leq j_0+10^{-4}$ then, $(S_n)_{n\in \Z}$ doesn't contain any of these subwords or their transpose.
\begin{proof}
By using the above boundaries (Lemma \ref{encadrement13}), it is simple computation.
    
\end{proof}
\begin{corollary}\label{cor1}
    If $S$ is a sequence such that $m(S) \leq j_0+10^{-4}$, then, the words $322,223$ and $323$ are forbidden.
\end{corollary}
\begin{proof}
    If $S$ contains $323$, then since according to the Lemma \ref{fw1_1}, the words $13, 232$ and $3323$ are forbidden, we can't extend the word $323$ to the left without making a forbidden word appear.
    
    If $S = \dots322\dots$, then according to the Lemma \ref{fw1_1}, since $13$ and $232$ are forbidden, the word $322$ must extend to the left in such way: $S = \dots3322\dots $
    
    However, since the words $13,23322$ and $33322$ are forbidden, we can't extend the word $3322$ without having a contradiction.
    The same reasoning gives that $223$ is forbidden.    
    
\end{proof}

\begin{lemma}\label{locuniq1}(Local uniqueness I):\\
Let $S \in \{1,2,3\}^{\Z}$ be such that $\lambda_0(S) \leq j_0 +10^{-4}$. Then $S$ or $S^T$ must take on of these forms;
\begin{enumerate}
    \item $S = 1^*$ or $2^*$ and $\lambda_0(S) <  j_0-10^{-2}$,
    \item $S = 33^*3$ and $\lambda_0(S) <  j_0-10^{-2}$,
    \item $S=1233^*212$ and $\lambda_0(S) <  j_0-10^{-3}$,
    \item $S=2333^*212$ or $3333^*212$ and $\lambda_0(S) <  j_0-10^{-3}$,
    \item $S=111233^*21112$  and $\lambda_0(S) <  j_0-10^{-4}$,
    \item $S=\omega_1^*$.
\end{enumerate}

\end{lemma}
As a consequence, if $S$ verifies $|\lambda_0(S)-j_0|\leq 10^{-4}$ and $m(S) \leq j_0 +10^{-4}$ then $S$ must be of the form: $$S=\omega_1^*.$$

\begin{proof}
Let $S$ verifies the conditions of the lemma.
If $S = 1^*,2^*$ then, we have $\lambda_0(S) <  j_0-10^{-2}$.
Otherwise, $S=3^*$ and since $13 $ and $232$ are forbidden words, then:
$$S = 33^*2 , 23^*3\:\: \text{or} \:\left\{ S = 33^*3\:\:\text{and}\:\: \lambda_0(S)<j_0-10^{-2}\right\}.$$ 
By symmetry, it is enough to study the case $S = 33^*2$. Since, $323$ and $322$ are forbidden, then $S= 33^*21$.

Again, $13$ is forbidden so:
$$S = 333^*21 \:\text{or}\: S=233^*21.$$
If $S = 333^*21$,
using the forbidden words $13 $ and $333211$, we must have:  $$S = 2333^*212 \:\:\text{or}\:\:S= 3333^*212 \:\:\text{and}\:\: \lambda_0(S) <  j_0-10^{-3}.$$

If  $S = 233^*21$,
using the forbidden words $323, 223$ and  $13$, we must have:
$$S = 1233^*211 \: \text{or} \:\left\{ S = 1233^*212\:\:\text{and}\:\: \lambda_0(S)<j_0-10^{-3}\right\}.$$

If $S = 1233^*211$, since $13, 31$ and $12332112$ are forbidden:

$$S = 11233^*2111 \:\text{or} \: S = 21233^*2111.$$

If $S = 11233^*2111$, since the words $13, 31, 21123321$ and $11123321111$ are forbidden, then:
$$ S = 111233^*21112\:\:\text{and}\:\: \lambda_0(S) <  j_0-10^{-4}.$$
\end{proof}

\begin{lemma}\label{fw2_1}(Forbidden words II):
\begin{enumerate}
    \item If $S = \omega_1^*\,1$, then $\lambda_0(S)> j_0+10^{-3}$.
    \item If $S = \,\omega_1^*\,22$ or \:$\,\omega_1^*\,23$, then $\lambda_0(S)> j_0+10^{-5}$.
    \item If $S = 21\,\omega_1^*\,211$ or \:$21\,\omega_1^*\,212$, then $\lambda_0(S)> j_0+10^{-6}$.
    \item If $S = 11\,\omega_1^*\,211$, then $\lambda_0(S)> j_0+10^{-5}$.
    \item If $S = 12111\,\omega_1^*\,2123$ or $22111\,\omega_1^*\,2123$, then $\lambda_0(S)> j_0+10^{-7}$.
    \item If $S = 332111\,\omega_1^*\,212333$ , then $\lambda_0(S)> j_0+10^{-8}$.
    \item If $S = 2332111\,\omega_1^*\,21233212$, then $\lambda_0(S)> j_0+8\times10^{-9}$.
    \item If $S = 112332111\,\omega_1^*\,2123321112$ , then $\lambda_0(S)> j_0+3\times10^{-10}$.
\end{enumerate}
    
\end{lemma}

\begin{corollary}
    If $S = 21\,\omega_1^*\,21$, then $m(S) > j_0 + 10^{-6}$.
    
    Therefore, if $m(S) <j_0+10^{-6}$, then the word $21\,\omega_1\,21$ is forbidden.
\end{corollary}
\begin{proof}
    Let $S = 21\,\omega_1^*\,21$. Then according to the Lemma \ref{fw2_1}, if $S_{7} \in \{1,2\}$, then $\lambda_0(S)>j_0+10^{-6}$.
    
    Otherwise, $S_{7} = 3$ and $\lambda_7(S) > j_0 +10^{-2}$. In both cases, we have $m(S) >j_0+10^{-6}$.
\end{proof}

\begin{lemma}\label{locuniq2}(Local uniqueness II)\\
Let $S \in \{1,2,3\}^{\Z}$ be such that
$$S = 21233^*2111 = \omega_1^*,$$
 and that $$\fa n \in \Z, \: \lambda_n(S) \leq j_0 +3\times10^{-10},$$
 then $S$ must take one of these forms;
\begin{enumerate}
    \item $S=33\,\omega_1^*\,21$ and $\lambda_0(S) <  j_0-10^{-4}$,
    \item $S=12\,\omega_1^*\,21$ and $\lambda_0(S) <  j_0-10^{-5}$,
    \item $S=22\,\omega_1^*\,212$  and $\lambda_0(S) <  j_0-10^{-6}$,
    \item $S=222\,\omega_1^*\,211,\, 211\,\omega_1^*\,212\:\text{or}\:\: 332\,\omega_1^*\,212$ and $\lambda_0(S) <  j_0-10^{-6}$,
    \item $S=1111\,\omega_1^*\,212$  and $\lambda_0(S) <  j_0-10^{-6}$,
    \item $S=\:1122\,\omega_1^*\,2111\:,\:1122\,\omega_1^*\,2112\:, \:2122\,\omega_1^*\,2111\:,\:2122\,\omega_1^*\,2112\:,\:2111\,\omega_1^*\,2121\:,\:2111\,\omega_1^*\,2122\:,$\\ 
    $2332\,\omega_1^*\,\,2111, 2332\,\omega_1^*\,2112, 3332\,\omega_1^*\,2111, 3332\,\omega_1^*\,2112\:$  and $\lambda_0(S) <  j_0-10^{-6}$,
    \item $S=2111\,\omega_1^*\,2123$.
\end{enumerate}
        
\end{lemma}
As a consequence, if $S$ verifies $|\lambda_0(S)-j_0|\leq 10^{-6}$ and $m(S) \leq j_0 +10^{-6}$ then $S$ must be of the form: $$S=2111\,\omega_1^*\,2123.$$

\begin{proof}
    Assuming that $S = 21233^*2111 = \omega_1^*$
    then since $\,\omega_1\,1$ and $13$ are forbidden words, $S$ must write itself in such way:
    $$S = 1\,\omega_1^*\,2, \: 2\,\omega_1^*2 \:\text{or}\: 3\,\omega_1^*2.$$
    
    \begin{enumerate}
        \item If $S= 3\,\omega_1^*2$, since the words $13, 232, \,\omega_1\,22$ and $\,\omega_1\,23$ are forbidden, then $S$ must write itself in such way:
    $$S = 33\,\omega_1^*21 \:\text{and}\: \lambda_0(S) <  j_0-10^{-4}.$$
        \item If $S = 2\,\omega_1^*\,2$, then, using the forbidden words $\omega_1\,22$ and $\omega_1\,23$ we have three possibilities:
    \begin{enumerate}
        \item $S = 22\,\omega_1^*\,21$,
        \item $S = 32\,\omega_1^*\,21$,
        \item $S = 12\,\omega_1^*\,21$ and $\lambda_0(S)<  j_0-10^{-5}$.
    \end{enumerate}
    We analyze each case separately:
    \begin{enumerate}
        \item If $S = 22\,\omega_1^*21$, then since $13$ and $322$ are forbidden we have:
        $$S = 122\,\omega_1^*\,211 \:\:\text{or}\:\:\left\{S =122\,\omega_1^*\,212, 222\,\omega_1^*\,212, 222\,\omega_1^*\,211 \:\:\text{and}\:\: \lambda_0(S)< j_0-10^{-6}\right\}.$$
        If $S = 122\,\omega_1^*\,211$, then since the words $13$ and $31$ are forbidden, we must have:
        $$S = 1122\,\omega_1^*\,2111,1122\,\omega_1^*\,2112,2122\,\omega_1^*\,2111,2122\,\omega_1^*\,2112 \:\:\text{and}\:\:\lambda_0(S)< j_0-10^{-5}.$$
        
        \item If $S = 32\,\omega_1^*\,21$ then since, $13, 31$ and $232$ are forbidden, we must have: 
        $$S =332\,\omega_1^*\,211 \:\text{or}\:\:\left\{ 332\,\omega_1^*\,212 \:\:\text{and}\:\:\  \lambda_0(S)<  j_0-10^{-4}\right\}.$$
        If $S =332\,\omega_1^*\,211$, then since $13, 31$ are forbidden, we must have:
        $$S = 2332\,\omega_2^*\,2111, 2332\,\omega_2^*\,2112, 3332\,\omega_2^*\,2111, 3332\,\omega_2^*\,2112 \:\:\text{and}\:\: \lambda_0(S)< j_0 -10^{-6}.$$
    \end{enumerate}
        \item If $S = 1\,\omega_1^*\,2$, since $31, \,\omega_1\,22$ and $\,\omega_1\,23$ are forbidden, so $S$ can take 2 forms:
    $$S = 11\,\omega_1^*\,21 \: \text{or} \: 21\,\omega_1^*\,21.$$

    However, since the words $21\,\omega_1\,211, 21\,\omega_1\,212$ and $13$ are forbidden, the word $21\,\omega_1\,21$ is not allowed.
    
    Hence, if $S = 11\,\omega_1^*\,21$, then, since the words $13$  and $11\,\omega_1\,211$ are forbidden, then :
    $$S = 111\,\omega_1^*\,212\: \text{or}\:\:\left\{  211\,\omega_1^*\,212 \:\: \text{and} \:\: \lambda_0(S)< j_0-10^{-5} \right\}.$$
    In the first case, since $13$ is forbidden,
    then:
    \begin{enumerate}
        \item $S = 1111\,\omega_1^*\,2121, 1111\,\omega_1^*\,2122 \: \text{or} \: 2111\,\omega_1^*\,2121$ and $\lambda_0(S)< j_0-10^{-5}$,
        \item $S = 2111\,\omega_1^*\,2122 \:\text{or} \: 1111\,\omega_1^*\,2123$ and $\lambda_0(S)< j_0-1.8\times10^{-6}$,
        \item $S = 2111\,\omega_1^*\,2123$.
    \end{enumerate}
    \end{enumerate}

\end{proof}

\begin{lemma}\label{fw3}(Forbidden words III (and last)):
\begin{enumerate}
    \item If $S = 3\,\omega_1\,\omega_1^*\,\omega_1\,21$, then $\lambda_0(S)> j_0+1.22\times10^{-10}$.
    \item If $S = 12\,\omega_1\,\omega_1^*\,\omega_1\,21$, then $\lambda_0(S)> j_0+2.7\times10^{-11}$.
    \item If $S = 22\,\omega_1\,\omega_1^*\,\omega_1\,212$, then $\lambda_0(S)> j_0+10^{-10}$.
    \item If $S = 222\,\omega_1\,\omega_1^*\,\omega_1\,21$, then $\lambda_0(S)> j_0+1.25\times10^{-12}$.
    \item If $S = 211\,\omega_1\,\omega_1^*\,\omega_1\,212$, then $\lambda_0(S)> j_0+3.15\times10^{-11}$.
    \item If $S = 1111\,\omega_1\,\omega_1^*\,\omega_1\,212$, then $\lambda_0(S)> j_0+1.9\times10^{-12}$.
    \item If $S = 111\,\omega_1\,\omega_1^*\,\omega_1\,2121$, then $\lambda_0(S)> j_0+2\times10^{-11}$.
    \item If $S = 2111\,\omega_1\,\omega_1^*\,\omega_1\,2122$, then $\lambda_0(S)> j_0+2.3\times10^{-12}$.
\end{enumerate}
\end{lemma}
   
\subsection{Sequence development Tree around $j_0 = \lambda_0(\overline{21233^*2111})$}
The previous work, of finding forbidden words and small words (i.e. with value less than $j_0$), can be represented visually by a tree. Indeed, assume that we are looking around the possible extension of a finite word $w_i\in \{1,2,3\}^{(\N)}$. In addition, suppose we have a set of forbidden words $F(w_i)$, which are all the forbidden words gathered by developing the sequence until the word $w_i$. Then, $\fa (s_l,s_r) \in \{1,2,3\}^2$ we have 4 possibilities:
\begin{enumerate}
    \item The word $s_l\,w_i\,s_r$ contains a forbidden word from $F(w_i)$. This case is immediately deleted.
    \item The word $s_l\,w_i\,s_r$ verifies $\min\lambda_0(\dots s_l\,w_i\,s_r\dots ) > j_0$, according to the bounding rules from Lemma \ref{encadrement13}. Then, $s_l\,w_i\,s_r$ is added to the new list of forbidden words. The branch associated with $s_l\,w_i\,s_r$ is then colored in red and ended.
    \item The word $s_l\,w_i\,s_r$ verifies $\max\lambda_0(\dots s_l\,w_i\,s_r\dots ) < j_0$ according to the bounding rules from Lemma \ref{encadrement13}. Then, the branch associated with $s_l\,w_i\,s_r$ is colored in blue and ended.
    \item The word $s_l\,w_i\,s_r$ verifies $\min\lambda_0(\dots s_l\,w_i\,s_r\dots ) < j_0 <\max\lambda_0(\dots s_l\,w_i\,s_r\dots )$ according to the bounding rules from Lemma \ref{encadrement13}. Then, the branch associated with $s_l\,w_i\,s_r$ is colored in black and must be extended in the next level.
    
\end{enumerate}
Using these rules, we construct a finite tree from the development of the sequence $\overline{212332111}$.

\subsection{Self-replication}
  We have gathered enough forbidden words to establish the self replication property of the function $\lambda_0$.
Let us call the whole set of forbidden words $F_{tot}$. We have:
$$F = \{13,232,323,322,333211,12332112,11123321111,\,\omega_1\,1,
 \,\omega_1\,23 ,\,\omega_1\,22,21\,\omega_1\,21,11\,\omega_1\,211, $$
 $$
 22111\,\omega_1\,2123,
 12111\,\omega_1\,2123,
 332111\,\omega_1\,212333,
 2332111\,\omega_1\,21233212,
 112332111\,\omega_1\,\,\omega_1\,2,
 3\,\omega_1\,\omega_1\,\omega_1\,21,$$
 $$ 12\,\omega_1\,\omega_1\,\omega_1\,21,
 222\,\omega_1\,\omega_1\,\omega_1\,21,
 22\,\omega_1\,\omega_1\,\omega_1\,212,
 211\,\omega_1\,\omega_1\,\omega_1\,212,
 1111\,\omega_1\,\omega_1\,\omega_1\,212,
 $$
 $$
 111\,\omega_1\,\omega_1\,\omega_1\,2121,
 2111\,\omega_1\,\omega_1\,\omega_1\,2122\},$$
and
$$F_{tot} = F\cup F^T.$$

If we use all of them, we will have self replication on both sides and end by showing that $j_0$ is isolated in $M$. To obtain a more comprehensive description of $M\cap (j_0,j_1)$, we need to consider the forbidden words with the smallest $\lambda_0$ value.

\subsubsection{Focus on forbidden words with the smallest $\lambda_0$ value}

We now focus on words from $F_{tot}$ with small $\lambda_0$ value to find the borders of $L$ and $M$. We can start by eliminating one forbidden word:
\begin{lemma}
    We have:
    \item $\sup_{n \in \Z} \lambda_n(\dots 222\,\omega_1\, \omega_1^*\,\omega_1\,21\dots ) > j_0 +10^{-10}$.
    
\end{lemma}
\begin{proof}
 We have:
 $$m(\dots 222\,\omega_1\, \omega_1^*\,\omega_1\,21\dots ) \geq \lambda_0(\dots 222\,\omega_1\, \omega_1^*\,\omega_1\,212\dots )>j_0 +10^{-10}.$$
\end{proof}

Let us consider 2 following words:
\begin{enumerate}
    \item $w_j = 2111\,\omega_1\, \omega_1\,\omega_1\,2122 \in F_{tot}$.
    \item $w_J = 1111\,\omega_1\, \omega_1\,\omega_1\,212 \in F_{tot}$.
\end{enumerate}
We also denote $w_j^*$ and $w_J^*$ respectively the words $2111\,\omega_1\, \omega_1^*\,\omega_1\,2122$ and $1111\,\omega_1\, \omega_1^*\,\omega_1\,212$  with the position zero marked with asterisk.
\begin{theorem}\label{sr_1}(Self-replication)

    Let $S  \in \{1,2,3\}^{\Z}$ such that: $S_{n-8}\dots S_{n+8} = 2111\,\omega_1^*2123$ for some $n \in \Z$.
    \begin{enumerate}
        \item If the finite word $S_{n-17}\dots S_{n+17}$ doesn't contain any words of $F_{tot}\setminus\{w_j,w_j^T\}$, then $S_{n-17}\dots S_{n+16} = 2111\,\omega_1\,\omega_1^*\,\omega_1\, 212$.
        \item If the finite word $S_{n-17}\dots S_{n+17}$ doesn't contain any word of $F_{tot}$,
    then $S_{n-17}\dots S_{n+17} = 2111\,\omega_1\,\omega_1^*\,\omega_1\, 2123$.
    \end{enumerate}

\end{theorem}

\begin{proof}
    Let assume that the finite word $S_{n-17}\dots S_{n+17}$ doesn't contain any word of $F_{tot}\setminus\{w_j,w_j^T\}$, and that we have:
    $$S_{-8}\dots S_{8} = 2111\,\omega_1^*\, 2123.$$
    Forbidden words used $12111\,\omega_1\,2123, 22111\,\omega_1\,2123$:
    $$S_{-9}\dots S_{8} = 32111\,\omega_1^*\, 2123.$$
    Forbidden words used $13 , 232$:
    $$S_{-10}\dots S_{8} = 332111\,\omega_1^*\, 2123.$$
    Forbidden words used $13 , 333211$:
    $$S_{-11}\dots S_{8} = 2332111\,\omega_1^*\, 2123.$$
    Forbidden words used $223,323$:
    $$S_{-12}\dots S_{8} = 12332111\,\omega_1^*\, 2123.$$
    Forbidden words used $31, 232$:
    $$S_{-12}\dots S_{9} = 12332111\,\omega_1^*\, 21233.$$
    Forbidden words used $31, 332111\,\omega_1\,212333$:
    $$S_{-12}\dots S_{10} = 12332111\,\omega_1^*\, 212332.$$
    Forbidden words used $322,323$:
    $$S_{-12}\dots S_{11} = 12332111\,\omega_1^*\, 2123321.$$
    Forbidden words used $13, 2332111\,\omega_1\,21233212$:
    $$S_{-12}\dots S_{12} = 12332111\,\omega_1^*\, 21233211.$$
    Forbidden words used $13, 12332112$:
    $$S_{-12}\dots S_{13} = 12332111\,\omega_1^*\, \omega_1.$$
    Forbidden words used $13, \omega_1\,1$:
    $$S_{-12}\dots S_{14} = 12332111\,\omega_1^*\, \omega_1\,2.$$
    Forbidden words used $\omega_1\,22, \omega_1\,23$:
    $$S_{-12}\dots S_{15} = 12332111\,\omega_1^*\, \omega_1\,21.$$
    Forbidden words used $11\,\omega_1\,211,13$:
    $$S_{-12}\dots S_{16} = 12332111\,\omega_1^*\, \omega_1\,212.$$
    Forbidden words used $31,112332111\,\omega_1\,\omega_1\,2$:
    $$S_{-13}\dots S_{16} = \omega_1\,\omega_1^*\, \omega_1\,212.$$
    Forbidden words used $3\,\omega_1\,\omega_1\,\omega_1\,21, 12\,\omega_1\,\omega_1\,\omega_1\,21, 22\,\omega_1\,\omega_1\,\omega_1\,212, 322$:
    $$S_{-14}\dots S_{16} = 1\,\omega_1\,\omega_1^*\, \omega_1\,212.$$
    Forbidden words used $31, 21\,\omega_1\,21$:
    $$S_{-15}\dots S_{16} = 11\,\omega_1\,\omega_1^*\, \omega_1\,212.$$
    Forbidden words used $31, 211\,\omega_1\,\omega_1\,\omega_1\,212 $:
    $$S_{-16}\dots S_{16} = 111\,\omega_1\,\omega_1^*\, \omega_1\,212.$$
    Forbidden words used $31, 1111\,\omega_1\,\omega_1\,\omega_1\,212 = w_J$:
    $$S_{-17}\dots S_{16} = 2111\,\omega_1\,\omega_1^*\, \omega_1\,212.$$
At this point we obtain the first result, since we only used words from the set $F_{tot}\setminus \{w_j\}$.

If $w_j$ is forbidden in the sequence $S_{n-17}\dots S_{n+17}$, then using the forbidden words, $ 111\,\omega_1\,\omega_1\,\omega_1\,2121, w_j = 2111\,\omega_1\,\omega_1\,\omega_1\,2122$ we must have:
$$S_{-17}\dots S_{17} = 2111\,\omega_1\,\omega_1^*\, \omega_1\,2123,$$ 
which proves the second point of the theorem.
\end{proof}

Using the Lemma \ref{sr_1} repeatedly, we have the two following results.
\begin{corollary}
    
\label{selfrep}
    Let $S  \in \{1,2,3\}^{\Z}$ be such that: $S_{-8}\dots S_{8} = 2111\,\omega_1^*2123$.
    \begin{enumerate}
        \item Then, if all the words from $F_{tot}\setminus\{w_j,w_j^T\}$ are forbidden in $S$, then we must have:
    $$S = \overline{\omega_1}\,\omega_1^*\,\omega_1\,212S_r \:\:\text{with}\:\: S_r\in \{1,2,3\}^{\N}.$$
        \item If all the words from $F_{tot}$ are forbidden in $S$, then we must have:
    $$S = \overline{\omega_1} = \overline{21233^*2111}.$$
    
    \end{enumerate}

\end{corollary}

\begin{lemma}\label{j_J}
We have:
\begin{enumerate}
    \item  
    \begin{align*}
    &\begin{multlined}
    \min\{ \lambda_0(S_l^T w_j^*S_r) / (S_l,S_r)\in \{1,2,3\}^{\N}\, \text{with}\:\lambda_0(S_l^T w_j^*S_r)=m(S_l^T w_j^*S_r)\:\:\text{and}\:\:S_l^T w_j^*S_r \:\:\text{doesn't contain}\\
    \text{any words from}\:\:
    F_{tot}\setminus\{w_j,w_j^T\}\}
    \end{multlined} \\
    &= \lambda_0(\overline{\omega_1}\,\omega_1\,\omega_1^*\,\omega_1\,2122\,\overline{12}) \in M.
    \end{align*}
    \item 
    \begin{align*}
    &\begin{multlined}
    \min\{\lambda_0(S_l^T \,w_J^* \, S_r)\: /\: (S_r,S_l) \in \{1,2,3\}^{\N}, \lambda_0(S_l^T \,w_J^* \, S_r) = m(S_l^T \,w_J^* \, S_r) \:\:\text{and}\:\:S_l^T w_J^*S_r \:\:\text{doesn't contain} \\
    \text{any words from}\:\:  F_{tot}\setminus\{w_j,w_j^T,w_J,w_J^T\}\}
    \end{multlined} \\
    &= \lambda_0(\overline{21}\,1111\,\omega_1\,\omega_1^*\,\omega_1\,\omega_1\,2122 \,1212\,  \overline{333212}) \in M.
    \end{align*}
    \end{enumerate}
    For the future, we call:
    $$j = \lambda_0(\overline{\omega_1}\,\omega_1\,\omega_1^*\,\omega_1\,2122\,\overline{12}) \:\: \text{and}\:\: J = \lambda_0(\overline{21}\,1111\,\omega_1\,\omega_1^*\,\omega_1\,\omega_1\,2122 \,1212\,  \overline{333212}).$$
    We get numerically:
    \begin{enumerate}
        \item $j \approx  j_0 + 8.32039\times 10^{-12}$.
        \item $J \approx j_0 + 8.42651\times 10^{-12}$. 
    \end{enumerate}

\begin{remark}
    We will show later that the number $j$ is the first element of the set $(M\setminus L) \cap(j_0,j_1)$ and that $J = j_1$ is the right border of the gap of $L$. 
\end{remark}

\end{lemma}
\begin{proof}
    \textbf{Proof of $j$}:
    
    Let $S = S_l^Tw_j^*S_r$ be such that $S$ doesn't contain any words from  $F_{tot}\setminus\{w_j,w_j^T\}$.
    Then, $S_{-8}\dots S_{8} = 2111\,\omega_1^*2123$ and according to the Corollary \ref{selfrep}, we must have $S = \overline{\omega_1}\,\omega_1^*\,\omega_1\,2122S_{18}\dots$.
    
    Since we are minimising the function $\lambda_0$, we must take $S_{18} = 1$ and because the word $13$ is forbidden, we must continue with the periodic sequence $\overline{12}$.
    The sequence that minimises $\lambda_0$ is:

    $$S= \overline{\omega_1}\,\omega_1\,\omega_1^*\,\omega_1\,2122\,\overline{12}
    =\overline{\omega_1}\,21233\,w_j^*\,\overline{12}.$$
    \textbf{Proof of $J$:}
    
    Let $S = S_l^Tw_J^*S_r$ such that $S$ doesn't contain any words from  $F_{tot}\setminus\{w_j,w_j^T,w_J,w_J^T\}$.
    Again, we have:
    $$S = \dots S_{-18}\,1111\,\omega_1\,\omega_1^*\,\omega_1\,212S_{17}\dots$$

    In order to minimise $\lambda_0$, since $13$ is forbidden, the left extension must be:
    $$\dots S_{18} = \overline{21}$$
    So:
    $$S = \overline{21}\,1111\,\omega_1\,\omega_1^*\,\omega_1\,212S_{17}\dots$$
    Again, in order to minimise $\lambda_0$, we must take $S_{17}= 3$ and therefore we have:
    $$S_{1}\dots S_{17} = 2111\,\omega_1\,2123$$
    Using the second point of the Theorem \ref{sr_1}, we must have:
    $$S = \overline{21}\,1111\,\omega_1\,\omega_1^*\,\omega_1\,\omega_1\,212S_{26}\dots$$
    Since the word $111\,\omega_1\,\omega_1\,\omega_1\,2121$ is forbidden (point 7 in Lemma \ref{fw3}) and that we minimise $\lambda_0(S)$, we must have $ S_{26} = 2$ and
    $$S = \overline{21}\,1111\,\omega_1\,\omega_1^*\,\omega_1\,\omega_1\,2122S_{27}\dots$$
    At this stage, we have a competition between the value of $\lambda_0(S)$ and $\lambda_9(S)$. 
    Indeed, we have to extend $S$ in such a way that $\lambda_9(S)$ is small enough so we can have $\lambda_9(S)\leq \lambda_0(S)$ and that in addition, $\lambda_0(S)$ is minimised. However, we know that:
    $$\lambda_9(\overline{21}\,1111\,\omega_1\,\omega_1^*\,\omega_1\,\omega_1\,2122 \,\overline{12}) = \lambda_0(\overline{21}\,1111\,\omega_1\,\omega_1\,\omega_1^*\,\omega_1\,2122 \,\overline{12}) <j.$$
    So we define, for $n\geq 0$:
    $$L_{n}^{(9)} = \max \left\{\lambda_9(\overline{21}\,1111\,\omega_1\,\omega_1^*\,\omega_1\,\omega_1\,2122 [12]_n T)\:\:/\:\: T\:\: \text{doesn't contain any words from} \:F_{tot}\setminus\{w_j,w_j^T,w_J,w_J^T\}\right\},$$
    $$l^{(0)}_n = \min \left\{\lambda_0(\overline{21}\,1111\,\omega_1\,\omega_1^*\,\omega_1\,\omega_1\,2122 [12]_n T)\:\:/\:\: T\:\: \text{doesn't contain any words from} \:F_{tot}\setminus\{w_j,w_j^T,w_J,w_J^T\}\right\},$$
    where $[12]_n$ represents the $n^{th}$ first terms of the sequence $(\overline{12})_{k\geq 0}$. 
    
Let us determinate the exact value of $l^{(0)}_{2n}$. We minimise the $\lambda_0$ value of sequence of following form $S = \overline{21}\,1111\,\omega_1\,\omega_1^*\,\omega_1\,\omega_1\,2122 (12)^nS_{27+2n}\dots$, so we must have:
$$S = \overline{21}\,1111\,\omega_1\,\omega_1^*\,\omega_1\,\omega_1\,2122 \,(12)^n 3S_{28+2n}\dots$$
But since the words $31$ and $232$ are forbidden, we must have:
$$S = \overline{21}\,1111\,\omega_1\,\omega_1^*\,\omega_1\,\omega_1\,2122 \,(12)^n 33S_{29+2n}\dots$$
Since the word $31$ is forbidden, to minimise $\lambda_0$, we must have:
$$S = \overline{21}\,1111\,\omega_1\,\omega_1^*\,\omega_1\,\omega_1\,2122\, (12)^n 333S_{30+2n}\dots$$
Since, the word $31$ is forbidden, and that we are minimising $\lambda_0$, we must have:
$$S = \overline{21}\,1111\,\omega_1\,\omega_1^*\,\omega_1\,\omega_1\,2122\, (12)^n 3332S_{31+2n}\dots$$
Since, the words $323$ and $322$ are forbidden, we must have:
$$S = \overline{21}\,1111\,\omega_1\,\omega_1^*\,\omega_1\,\omega_1\,2122\, (12)^n 33321S_{32+2n}\dots$$
Since the word $333211$ is forbidden, we must have:
$$S = \overline{21}\,1111\,\omega_1\,\omega_1^*\,\omega_1\,\omega_1\,2122\, (12)^n 333212S_{33+2n}\dots$$

Then, we repeat the process by induction, since $2n +33$ and $2n+27$ have the same parity. So finally, we have:
$$S = \overline{21}\,1111\,\omega_1\,\omega_1^*\,\omega_1\,\omega_1\,2122\, (12)^n \, \overline{333212}.$$
And:
$$l_{2n}^{(0)} = \lambda_0\left(\overline{21}\,1111\,\omega_1\,\omega_1^*\,\omega_1\,\omega_1\,2122\, (12)^n \, \overline{333212}\right).$$

Now we determine the exact value of $l^{(0)}_{2n+1}$. If we minimise the $\lambda_0$ value of the sequence with the form $S = \overline{21}\,1111\,\omega_1\,\omega_1^*\,\omega_1\,\omega_1\,2122\, (12)^n\,1\,S_{2n+28}\dots $, then we must have:
$$S = \overline{21}\,1111\,\omega_1\,\omega_1^*\,\omega_1\,\omega_1\,2122\, (12)^n 1\,1\dots$$
And then, using the fact that the words $13,31$ are forbidden, we must continue with the periodic sequence $\overline{12}$, so we have:
$$l_{2n+1}^{(0)} = \lambda_0\left(\overline{21}\,1111\,\omega_1\,\omega_1^*\,\omega_1\,\omega_1\,2122\, (12)^n 1\, \overline{12}\right).$$

 Since $0$ is even and $9$ is odd, we have that 
    $\fa n\in \N:$
    $$L_{2n}^{(9)} = \lambda_9\left(\overline{21}\,1111\,\omega_1\,\omega_1^*\,\omega_1\,\omega_1\,2122\, (12)^n \, \overline{333212}\right) \:\:\text{and}\:\: L_{2n+1}^{(9)} = \lambda_0\left(\overline{21}\,1111\,\omega_1\,\omega_1^*\,\omega_1\,\omega_1\,2122\, (12)^n \, 1\,\overline{12}\right).$$
    
    So finally, we have that:
    $$J = \min_{n \in \N}\left\{\max\left(l_n^{(0)}, L_n^{(9)}\right)\right\}.$$
    The sequence $(L_n^{(9)})_{n\geq 0}$ is decreasing to $L^{(9)}$ and the sequence $(l_n^{(0)})_{n\geq 0}$ is increasing to $l^{(0)}$ with:
    $$L^{(9)}<j<l^{(0)}.$$
    If we define  $N = \min\{n\in \N \:/\: l_n^{(0)}\geq L_n^{(9)}\}$, then we have:
    $$J = l^{(0)}_N.$$
    We get computationally that $N = 4$ and therefore:
    $$J = l^{(0)}_4 = \lambda_0(\overline{21}\,1111\,\omega_1\,\omega_1^*\,\omega_1\,\omega_1\,2122 \,1212\,  \overline{333212}) \approx j_0 +  8.42651\times 10^{-12}.$$
\end{proof}

\begin{prop}\label{remark fw}
    Let $S\in \{1,2,3\}^{\Z}$ be a sequence such that $\lambda_0(S) = m(S) = \sup_{n \in \Z} \lambda_n(S)$.
    
    According to the definition of $j$ and $J$ in Lemma \ref{j_J}, we have:

\begin{enumerate}
    \item If $\lambda_0(S) \in M\cap (j_0,J+\alpha)$, for $\alpha >0$ small enough, then at least one of the words $w_j,w_j^T$ and $w_J,w_J^T$ must appear in $S$.
    \item  If $\lambda_0(S) \in M\cap(j_0,J)$, no forbidden words from $F_{tot}\setminus \{w_j,w_j^T\}$ are allowed to appear in $S$.
    \item If $\lambda_0(S) \in M\cap(j_0,j)$, no forbidden words from $F_{tot}$ are allowed to appear in $S$.
    
\end{enumerate}

\end{prop}

\begin{proof}
    Let $m\in M$ such that $m>j_0$ and $S\in \{1,2,3\}^{\Z}$ a sequence such that: $$m=\lambda_0(S) = \sup_{n\in \Z}\lambda_n(S).$$
    If $m>J$, then clearly one of $w_j,w_j^{T}$ and $w_J,w_J^{T}$ must appear in the sequence $S$.

    Let assume that $m \in (j_0,J)$. By contradiction, suppose that there exists $w_f\in F_{tot}\setminus \{w_j,w_j^T\}$ such that $S =  S_l^Tw_fS_r$. If $w_f \in F_{tot}\setminus\{w_j,w_j^T,w_J,w_J^T\}$, then according to Lemmas \ref{fw1_1}, \ref{fw2_1} and \ref{fw3}, we must have $\lambda_N(S) >j_0+2\times 10^{-11}$ for $N \in \Z$, which is impossible. \textbf{So $S$ doesn't contain any word from $F_{tot}\setminus\{w_j,w_j^T,w_J,w_J^T\}$}. 
    
    If $S$ contains $w_J$, then we have $S = S_l^Tw_JS_r$. By definition of $J$ (cf. Lemma \ref{j_J}), we must have $\lambda_N(S) \geq J$ for $N \in \Z$ and therefore, $m\geq J$, which is impossible. The same is true with $w_J^T$.
    \textbf{So $S$ doesn't contain any subwords from $F_{tot}\setminus\{w_j,w_j^T\}$.}
    
    Let assume that $m \in (j_0,j)$. Then according to above, $S$ doesn't contain any words from $F_{tot}\setminus\{w_j,w_j^T\}$.
    If $S = S_l^T\,w_j\,S_r$ then by definition of $j$ (cf. Lemma \ref{j_J}), we must have $\lambda_N(S)\geq j$ and therefore: $m(S) \geq j$. Again, the same is true for $w_j^T$. Hence, the sequence $S$ doesn't contain any words from $F_{tot}$. 
\end{proof}

The main consequence of all of this is the following Theorem:
\begin{theorem}
    
    We have:
    \begin{enumerate}
        \item $M\cap(j_0,j) = \emptyset$.
        \item $L\cap(j_0,J) = \emptyset$.
    \end{enumerate}
    
\end{theorem}

\begin{proof}
    Let assume by contradiction that there exist $m\in M\cap(j_0,j)$. Then we can find a sequence $S\in \{1,2,3\}^{\Z}$ such that $m = \sup_{n\in \Z}\lambda_n(S) = m(S) = \lambda_0(S)$.\:
    Then, we have $j_0<\lambda_0(S) < j< j_0 +  3\times 10^{-10}$.
    
    So according to the Lemma \ref{locuniq2}, we must have:
    $S_{-8}\dots S_{8} = 2111\,\omega_1^*2123$.
    
    In addition, according to the Proposition \ref{remark fw}  the words from the set $F_{tot}$ are forbidden in $S$. 
    Then, because of the Corollary \ref{selfrep}, the sequence $S$ must extend in such way: 
    $$S = \overline{\omega_1},$$ 
    and so $m(S) = j_0$, a contradiction.
    
    Let assume by contradiction that there exists $l\in L \cap (j_0,J)$. We use the fact that the Markov values of periodic sequences is dense in $L$ (see \cite[Theorem 2, Chapter 3]{Cusick-Flahive}). Therefore, $\exists \,(l_n)_{n\in \N} \in L^{\N}$ such that: 
    $$\fa n \in \N, \exists \,\sigma^{(n)} \in \N^{(\N)}, l_n = m\left(\overline{\sigma^{(n)}}\right).$$ And:
    $$\lim_{n\rightarrow \infty} l_n = l.$$
    Let us write $\fa n \in \N, S^{(n)} = \overline{\sigma^{(n)}}$. We can assume (even if it means taking $n \geq n_0$ with $n_0$ big enough), $\fa n \in \N, \,S^{(n)} \in \{1,2,3\}^{\Z}$, $l_n<J$.
    So we have $\lambda_0(S^{(n)}) = l_n <J$. Using the Lemma \ref{locuniq1} and \ref{locuniq2}, we have that $S_{-8}^{(n)}\dots S_8^{(n)} = 2111\,\omega_1\,2123$ (or its transpose).
    Since $l_n = \lambda_0(S^{(n)}) \in M\cap(j_0,J)$, according to the Proposition \ref{remark fw}, no forbidden words from $F_{tot}\setminus \{w_j,w_j^T\}$ are allowed to appear in $S^{(n)}$.
    Therefore, using the Corollary \ref{selfrep}, we have:
    $$\fa n \in \N \:\: S^{(n)} = \overline{\omega_1}\,\omega_1^*\,\omega_1\,212\,S_r^{(n)} \:\:\text{with}\:\:S_r^{(n)} \in \{1,2,3\}^{\N}.$$
    Hence, $\fa n \in \N,\:\: \overline{\sigma^{(n)}} = \overline{\omega_1}\,\omega_1^*\,\omega_1\,212\,S_r^{(n)}$ and necessarily $\fa n \in \N,\:\: \overline{\sigma^{(n)}} = S^{(n)} = \overline{\omega_1}$. So, for all $n\in\N$, we have $l_n = m(\overline{\omega_1}) = j_0$ and $l = j_0$, which is impossible.
\end{proof}

\subsection{The description of $(M\setminus L)\cap(j_0,j_1)$}

Now, let us characterize the set $(M\setminus L)\cap(j_0,J)$.
\begin{lemma}\label{extension}
We have:
    \begin{enumerate}
        \item $\lambda_0(\overline{\omega_1}\,\omega_1^*\,\omega_1\, 21223\dots)>\lambda_0(\overline{\omega_1}\,\omega_1^*\,\omega_1\, 21222\dots)> J$,
        \item $\lambda_0(\overline{\omega_1}\,\omega_1^*\,\omega_1\, 212211\dots)> J$,
        \item $\lambda_0(\overline{\omega_1}\,\omega_1^*\,\omega_1\, 2122123\dots)>\lambda_0(\overline{\omega_1}\,\omega_1^*\,\omega_1\, 2122122\dots)> J$,
        \item  $\lambda_0(\overline{\omega_1}\,\omega_1^*\,\omega_1\, 21221211\dots)> J$.
    \end{enumerate}
\end{lemma}

\begin{corollary}\label{carac1}
    Let $m \in M\cap(j,J)$ and $S \in \{1,2,3\}^{\Z}$ be a sequence such that $m = \lambda_0(S)$. Then we have:
    $$S = \overline{\omega_1}\,\omega_1^*\,\omega_1\, 21221212\, S_r,$$
    is such that:
    \begin{enumerate}
        \item $1212S_r$ doesn't contain any words from $F_{tot}\setminus\{w_j,w_j^T\} \cup \{2111\,\omega_1\,2123\}$, 
        \item If $\exists \,n \in\Z$ such that $S_{n-8}\dots S_{n+8} = 3212\,\omega_1^T\,1112$, then we have $S = \dots S_{n-16}\,212\,\omega_1^T\,\overline{\omega_1^T}$.
     \end{enumerate}
 \end{corollary}
\begin{proof}
    Let $m(S) = \lambda_0(S) \in (j_0,J)$. According to the Lemma \ref{locuniq1} and  \ref{locuniq2}, we have:
    $$S_{-8}\dots S_8 = 2111\,\omega_1^*\,2123.$$
    Since $m(S) < J$, the sequence $S$ doesn't contain any words from $F_{tot}\setminus\{w_j,w_j^T\}$. Then, because of the Theorem \ref{selfrep}, we have:
    $$S = \overline{\omega_1}\,\omega_1^*\,\omega_1\,212\dots$$
    If $S = \overline{\omega_1}\,\omega_1^*\,\omega_1\,2121\dots$, then according to point 7 of Lemma \ref{fw3}, $\lambda_0(S)> j_0 +2\times 10^{-11} >J$, which is impossible. If $S = \overline{\omega_1}\,\omega_1^*\,\omega_1\,2123\dots$, then $\lambda_0(S) < \lambda_0(\overline{\omega_1}\,\omega_1^*\,\omega_1\,2122\overline{12})=j$
    which again, is impossible. 
    
    So $S = \overline{\omega_1}\,\omega_1^*\,\omega_1\,2122\dots$. However, since we have $ \lambda_0(S) < J$, using the forbidden word $13$ and the Lemma \ref{extension}, we must have:
    $$S = \overline{\omega_1}\,\omega_1^*\,\omega_1\,21221212\,S_r.$$
    \begin{itemize}
    \item If $1212S_r$ contains any words from $F_{tot}\setminus\{w_j,w_j^T\}$,  we have obviously $m(S)>J$.
    
    \item If $1212S_r$ contains $2111\,\omega_1\,2123$, then using self-replication according to the Theorem \ref{selfrep}, we have a contradiction in the writing of $S$ since $1212$ is not a sub-word of $\overline{\omega_1}$.

    \item If $S$ contains $3212\,\omega_1^T\,1112$, we have $n \in \Z$ such that $S_{n-8}\dots S_{n+8} = 3212\,\omega_1^T\,2111$ then because of self-replication (Theorem \ref{sr_1}), applied with $\omega_1^T$, we must have $S_{n-16}\dots S_{n+4}\dots  = 212\,\omega_1^T\,\overline{\omega_1^T}$.
    \end{itemize}
\end{proof}

We can now make a global description of the set $(M \setminus L)\cap(j,J)$.

First, let us call $\Tilde{F}$ the set of the forbidden words (from $F_{tot}$) that doesn't contain the subword $2111\,\omega_1\,2123$. We call $\Tilde{F}^T$ the set of forbidden words (from $F_{tot}^T$) that doesn't contain the subword $3212\,\omega_1^T\,1112$.
We have, according to Lemma \ref{fw1_1}, \ref{fw2_1} and \ref{fw3}:
$$\Tilde{F_0} = \left\{13 , 232 , 323 , 322 , 333211 , 12332112 , 11123321111, \omega_1\,1, \omega_1\,22, \omega_1\,23 , 21\,\omega_1\,21, 11\,\omega_1\,211 \right\}.$$
We can improve the set of forbidden words:
\begin{lemma}
    Let $S \in \{1,2,3\}^{\Z}$.
    \begin{enumerate}
        \item If $S = \dots3321111\dots$,  then, $m(S) >j_0 +10^{-4}$.
        \item If $S = \dots332112\dots$,  then, $m(S) >j_0 +10^{-4}$.
    \end{enumerate}
   
\end{lemma}
\begin{proof}
    Let $S\in \{1,2,3\}\in{\Z}$ with $S = \dots3321111\dots$.
    Then using the forbidden words $31$ and $333211$, we have, according to Lemma \ref{fw1_1}:
    $$m(S) >j_0 + 10^{-3}\:\:\text{or}\:\: S =\dots23321111\dots $$
    In the last case, using the forbidden words $323$ and $322$, according to the Lemma \ref{fw1_1} and \ref{fw2_1} we have:
    $$m(S) >j_0 + 10^{-4}\:\:\text{or}\:\: S =\dots123321111\dots $$
    In the last case, using the forbidden words $31$ and $\omega_1\,$, according to the Lemma \ref{fw1_1} and \ref{fw3} we have:
    $$m(S) >j_0 + 10^{-3}\:\:\text{or}\:\: S =\dots1123321111\dots $$
    In the last case, using the forbidden words $31, (12332112)^T$ and $11123321111$, according to the Lemma \ref{fw1_1} we have:
    $$m(S) >j_0 + 10^{-4}.$$
    Hence, in every case we have: $m(S)>j_0+10^{-4}$. Using the same reasoning with the forbidden words $13, 333211, 322, 323$ and  $12332112$ we see that:
    If $S = \dots332112$, then $m(S)>j_0+10^{-4}$.
\end{proof}

The new set of forbidden words is:
$$\Tilde{F_1} = \left\{13 , 232 , 323 , 322 , 333211 , 332112 , 3321111, \omega_1\,22, \omega_1\,23 , 21\,\omega_1\,2, 1\,\omega_1\,211 \right\}.$$
We finally add the self-replicating word $2111\,\omega_1\,2123$ and its transpose $3212\,\omega_1^T\,1112$, so we get:
$$\Tilde{F_2} = \Tilde{F_1}\cup \Tilde{F_1}^T\cup \{2111\,\omega_1\,2123, 3212\,\omega_1^T\,1112\}.$$

\begin{theorem}\label{description globale M/L}
    We have:
    $$(M\setminus L)\cap (j,J) = C\cup D\cup X,$$
    where
    $$X = \left\{\lambda_0(\overline{\omega_1}\,\omega_1^*\,\omega_1\,2122\,121\,2212\,\omega_1^T\,\overline{\omega_1^T}), \lambda_0(\overline{\omega_1}\,\omega_1^*\,\omega_1\,2122\,12121\,2212\,\omega_1^T\,\overline{\omega_1^T}), \lambda_0(\overline{\omega_1}\,\omega_1^*\,\omega_1\,2122\,1212121\,2212\,\omega_1^T\,\overline{\omega_1^T})\right\},$$
    $$D = \{\lambda_0(\overline{\omega_1}\,\omega_1^*\,\omega_1\,21221212\, s\, 21212212\,\omega_1^T\,\overline{\omega_1^T}) \:\:/\:\:212\,s\,212 \in \{1,2,3\}^{(\N)} \:\text{doesn't contain any words}$$
    $$\text{from } \Tilde{F_2} \:\text{and}\:\: [0,s] \leq [0,s^T]\},$$
    and
    $$C = \left\{\lambda_0(\overline{\omega_1}\,\omega_1^*\,\omega_1\,21221212\, S)\:\: / \:\:S \in \{1,2,3\}^{\N}\:\text{doesn't contain words from } \Tilde{F_2} \right\}.$$
    \end{theorem}
    \begin{proof}
        Let $S\in \{1,2,3\}^{\Z}$ be a sequence such that $m(S) = \lambda_0(S) \in M\cap(j,J)$. Then according to the Corollary \ref{carac1} we have $S = \overline{\omega_1}\,\omega_1^*\,\omega_1\,21221212\,S_r$, with $S_r \in \{1,2,3\}^{\N}$ not containing any words from $F_{tot}\setminus \{w_j,w_j^T\}\cup \{2111\,\omega_1\,2123\}$.

        If $S$ doesn't contain the word $ 3212\,\omega_1^T\,1112$, then the set of forbidden words can be simplified and the sequence $S$ doesn't contain any word from $\Tilde{F_2}$. So $m(S) \in C$. 
        
        Otherwise, let $ N=\min\{n \in \N\: /\: 3212\,\omega_1^T\,1112 \:\text{ appears in}\: (S_{k})_{k\geq n} \}$. Then, we have $S_{N}\dots S_{N+16} = 3212\,\omega_1^T1112$ and because of Corollary \ref{carac1}, we must have:
        $$S_{N-8}\dots S_{N}\dots  = 212\,\omega_1^T\,\overline{\omega_1^T}.$$
        
        By definition of $N$ we must have $S_{N-9} \in \{1,2\}$ and using the Lemma \ref{extension}, we must have:
        $$S_{N-13}\dots  =21212212\,\omega_1^T\,\overline{\omega_1^T}.$$
        Then we have:
        \begin{enumerate}
            \item $S = \overline{\omega_1}\,\omega_1^*\,\omega_1\,21221212212\,\omega_1^T\,\overline{\omega_1^T}$,\, $(N = 30)$,
            \item $S = \overline{\omega_1}\,\omega_1^*\,\omega_1\,2122121212212\,\omega_1^T\,\overline{\omega_1^T}$,\, $(N = 32)$,
            \item $S = \overline{\omega_1}\,\omega_1^*\,\omega_1\,212212121212212\,\omega_1^T\,\overline{\omega_1^T}$,\, $(N = 34)$,
            \item $S = \overline{\omega_1}\,\omega_1^*\,\omega_1\,21221212\,s\, 21212212\,\omega_1^T\,\overline{\omega_1^T}$,\, $(N \geq 35)$.
        \end{enumerate}
        In the last case, we also need to add the conditions:
        $$\lambda_0\left(\overline{\omega_1}\,\omega_1^*\,\omega_1\,21221212\,s\, 21212212\,\omega_1^T\,\overline{\omega_1^T}\right) \geq \lambda_0\left(\overline{\omega_1}\,21221212\,s\, 21212212\,\omega_1^T\,(\omega_1^T)^*\,\overline{\omega_1^T}\right),$$
        which is equivalent to $[0,s]\leq [0,s^T]$,
        and also that the sequence $212\,s\,212$ doesn't contain any words from the simplified set $\Tilde{F_2}$. So $m(S) \in D$ and $(M\setminus L)\cap(j,J) \subset C\cup D \cup X$.

        Now we prove the reverse inclusion. We start by assuming $m\in C$. Then $m = \lambda_0(S)$ with $S = \overline{\omega_1}\,\omega_1^*\,\omega_1\,21221212\, S_r$ and $1212S_r$ not containing any words from $\Tilde{F}\cup \Tilde{F}^T\cup \{2111\,\omega_1\,2123,3212\,\omega_1^T\,1112\}$. Clearly, we have $ m \in (j,J)$. We want to show that $m = \sup_{n\in \Z}\lambda_n(S)$ and therefore $m \in M \cap (j,J)$.

        Firstly, we have, $\fa k \in \N^*$, $\lambda_{-9k}(S) = \lambda_0(\dots\,\omega_1\,\omega_1\,\omega_1^*\,\omega_1\,\omega_1\,\dots)< j_0 +10^{-14} < j$.
        
        We also have $\lambda_9(S) = \lambda_0(\dots2111\,\omega_1\,2122\dots) <j_0-10^{-6}$ according to the Lemma \ref{locuniq2}, so:
        $$\sup_{n\leq 20} \lambda_n(S) \leq \lambda_0(S).$$
        If we assume by contradiction that $\exists \,n \geq 21$ such that $\lambda_n(S) > \lambda_0(S)>j$, since $S$ doesn't contain any words from $\Tilde{F_2}$ and that $\lambda_n(S) >j_0$, using Lemmas \ref{locuniq1} and \ref{locuniq2} (forbidding the cases where $\lambda_0(S)<j_0$), we see that $S$ must take one of these forms:
        $$S_{n-8}\dots S_{n+8} = 2111\,\omega_1^*\,2123 \:\:\text{or}\:\: 3212\,(\omega_1^*)^T\,1112,$$
        which is impossible by definition of $S$. Then, $\fa n \geq 21, \lambda_n(S) \leq \lambda_0(S)$ and therefore $\lambda_0(S) = m(S)$. So $C\subset (M\setminus L)\cap(j,J)$. The same reasoning gives $D \subset (M\setminus L)\cap(j,J)$.

    \end{proof}

\begin{lemma}\label{isolated1}
    The set $D$ previously defined is a set of isolated points of $M\setminus L$.
\end{lemma}
\begin{proof}
    Let $m\in D$ and $s$ be a finite sequence such that $m = \lambda_0(\overline{\omega_1}\,\omega_1^*\,\omega_1\,21221212\, s\, 21212212\,\omega_1^T\,\overline{\omega_1^T})$, with $s$ such that $212\,s\,212 \in \{1,2,3\}^{(\N)}$ doesn't contain any words from  $\Tilde{F}\cup \Tilde{F}^T\cup \{2111\,\omega_1\,2123,3212\,\omega_1^T\,1112\}$.
    
    Let assume there is a sequence $(m_n)_{n \in\N}\in M^{\N}$ such that:
    $$m = \lim_{n \rightarrow \infty}m_n.$$
    We can assume that we have in addition $\fa n \in \N, m_n\in M\cap (j,J)$.
    Hence, according to the Corollary \ref{carac1}, $\fa n \in \N, \exists \,S^{(n)}$ a sequence such that: 
    $$m_n = \lambda_0(\overline{\omega_1}\,\omega_1^*\,\omega_1\,21221212\,S^{(n)}).$$
    Using the Proposition \ref{continuite lambda},  we can find a integer $n_0$ such that $\fa n\geq n_0$, we have: 
    $$S^{(n)} = s\,212\,\omega_1^T\omega_1^T\omega_1^T\, \Tilde{S}^{(n)}.$$
    So the sequence $S^{(n)}$ contains the words $3212\,\omega_1^T1112$ and according to the Corollary \ref{carac1}, we must have:
    $$\fa n \geq n_0 \:\: S^{(n)} = s\,21212212\,\overline{\omega_1^T}.$$
    Therefore $\fa n \geq n_0, \:\: m_n = m$. Hence, $m$ is an isolated point in $(M\setminus L) \cap (j,J)$.
\end{proof}
\begin{corollary}
    We have :
    
    $$\max (M \setminus L)\cap (j, J) = \lambda_0(\overline{\omega_1}\,\omega_1^*\,\omega_1\,212212\, \overline{123332})  \approx j + 4.4064196\times 10^{-14}.$$
    Thus:
    $$(M\setminus L)\cap (j_0,J) \subset (j, j + 4.4064196\times 10^{-14}).$$
\end{corollary}
\begin{proof}

This is a computation, we want to solve the following problem:
$$\max\left\{\lambda_0(\overline{\omega_1}\,\omega_1^*\,\omega_1\,21221212\, S_r) / S_r \:\text{doesn't contain forbidden words from}\: \Tilde{F_1} \cup \Tilde{F_1}^T\cup\{2111\,\omega_1\,2123\}\right\}.$$
Let $S = \overline{\omega_1}\,\omega_1^*\,\omega_1\,21221212\, S_{22}\dots$ be a sequence that maximise $\lambda_0$ and that doesn't contain forbidden words from $\Tilde{F_1} \cup \Tilde{F_1}^T\cup\{2111\,\omega_1\,2123\}$.
We have:
$$S = \overline{\omega_1}\,\omega_1^*\,\omega_1\,21221212\,3 S_{23}\dots$$
Then using the forbidden words $31, 232, 323, 322$ and $333211$, we must have:
$$S = \overline{\omega_1}\,\omega_1^*\,\omega_1\,21221212\,333212 S_{28}\dots$$
Using the same argument repeatedly, since both $22$ and $28$ are even, we have:
$$S = \overline{\omega_1}\,\omega_1^*\,\omega_1\,212212\, \overline{123332}.$$
\end{proof}

\subsection{The local border of $L$}

\begin{theorem}
    $J \in L'\subset L$ and therefore $(j_0,J)$ is the largest gap to the right of $j_0$ and we found that $j_1 = J$.
    \end{theorem}

    \begin{remark}
        With the notation of Theorem \ref{principal1}, this shows that $A_1 = (j_0,j_1)$.
    \end{remark}

\begin{proof}
By the proof of Lemma \ref{j_J} we have
\begin{align}\label{eq:J_char}
J &= 
\begin{multlined}[t]
\min \big\{\lambda_0(S_l^T \,1111\,\omega_1\omega_1^*\omega_1\omega_1\, 2122\, 1212 \, S_r)\: /\: (S_r,S_l) \in \{1,2,3\}^{\N}\:\text{and} \:S_l^T w_JS_r \:\:\text{doesn't contain} \\ \text{any words from}\:\: F_{tot}\setminus\{w_j,w_j^T,w_J,w_J^T\}\big\} 
\end{multlined}  \\
&= \lambda_0(\overline{21}\,1111\,\omega_1\,\omega_1^*\,\omega_1\,\omega_1\,2122 \,1212\,  \overline{333212}) \approx j_0 + 8.42651\times 10^{-12}. \nonumber  
\end{align}

    Let us call, for all $n\geq 2$ big enough:
    $$S^{(n)} = \overline{(21)^n\,1111\,\omega_1\,\omega_1^*\,\omega_1\,\omega_1\,2122\,1212\,(333212)^n}.$$

    We have:
    \begin{enumerate}
        \item $\lambda_{-9}(S^{(n)}) = \lambda_0(\dots1111\,\omega_1^*\,212\dots) <j_0-10^{-6}$ according to point 5 of Lemma \ref{locuniq2}.
        \item $\lambda_0(S^{(n)}) \geq J $ by the above characterization of $J$.
        \item $\lambda_9(S^{(n)}) = \lambda_0(\dots S_{-27-2n}(21)^n\,1111\,\omega_1\,\omega_1\,\omega_1^*\,\omega_1\,2122\,1212\,(333212)^n\,S_{22+6n}\dots)$$ \\
        \leq \lambda_0(\overline{21}\,1111\,\omega_1\,\omega_1\,\omega_1^*\,\omega_1\,2122\,1212\,(333212)^n\,3\dots) \leq \lambda_0(\overline{21}\,1111\,\omega_1\,\omega_1\,\omega_1^*\,\omega_1\,2122\,1212\,\overline{333212})+ \varepsilon(n)$.

        By direct computation we have:
        $\lambda_0(\overline{21}\,1111\,\omega_1\,\omega_1\,\omega_1^*\,\omega_1\,2122\,1212\,\overline{333212}) < J$.
        
        Therefore, for $n$ big enough, we have $\varepsilon(n) \leq J-\lambda_0(\overline{21}\,1111\,\omega_1\,\omega_1\,\omega_1^*\,\omega_1\,2122\,1212\,\overline{333212})$ and so $\lambda_9(S^{(n)}) < J$.
        \item $\lambda_{18}(S^{(n)}) = \lambda_0(\dots2111\,\omega_1^*\,2122\dots) <j_0-10^{-6}$ according to the point $6$ of Lemma \ref{locuniq2}.
        \item Lastly, $\fa \,0\leq k < n $: $$\max_{i\in \{-1,0,1\}}\lambda_{6k+32+ i}(S) = \lambda_{6k+31}(S) = \lambda_0(\dots2123^*33212\dots)<j_0 -10^{-3},$$
        according to point 4 from Lemma \ref{locuniq1}.
    \end{enumerate}
    So we have:
    $$l_n = m\left(\overline{(21)^n\,1111\,\omega_1\,\omega_1^*\,\omega_1\,\omega_1\,2122\,1212\,(333212)^n}\right) = l\left(\overline{(21)^n\,1111\,\omega_1\,\omega_1^*\,\omega_1\,\omega_1\,2122\,1212\,(333212)^n}\right) $$
    $$ = \lambda_0\left(\overline{(21)^n\,1111\,\omega_1\,\omega_1^*\,\omega_1\,\omega_1\,2122\,1212\,(333212)^n}\right). $$
    Therefore, $\fa n\geq 2$, we have $l_n \in L$,
    and we see that: 
    $$\lim_{n\rightarrow \infty}l_n = J.$$
    Since $L$ is closed, we have $J\in L$.
    
\end{proof}

\begin{theorem}
    We have, $\fa \alpha>0, \: HD(M\cap(j_1,j_1+\alpha)) =1$.
\end{theorem}
\begin{proof}

    Let $\alpha>0$. Using the continuity of the function $\lambda_0$, there is a $n_0 \in \N$ such that
    $\fa n>n_0$, $\fa S_l,S_r \in \{1,2\}^{\N}$:
    $$\lambda_0\left(S_l^T\,(21)^n\,1111\,\omega_1\,\omega_1^*\,\omega_1\,\omega_1\,2122\,1212(333212)^n\,S_r\right) \in [j_1,j_1+\alpha).$$
    In addition, we have:
    $$\lambda_0\left(S_l^T\,(21)^n\,1111\,\omega_1\,\omega_1^*\,\omega_1\,\omega_1\,2122\,1212(333212)^n\,S_r\right) = m\left(S_l^T\,(21)^n\,1111\,\omega_1\,\omega_1^*\,\omega_1\,\omega_1\,2122\,1212(333212)^n\,S_r\right)$$

    Indeed, if we call $S = S_l^T\,(21)^n\,1111\,\omega_1\,\omega_1^*\,\omega_1\,\omega_1\,2122\,1212(333212)^n\,S_r$, we have:
    \begin{enumerate}
        \item $\lambda_{-9}(S) = \lambda_0(\dots1111\,\omega_1^*\,212\dots) <j_0-10^{-6}$ according to point 5 of Lemma \ref{locuniq2}.
        \item $\lambda_0(S) \geq J$ from the characterization \eqref{eq:J_char}.
        \item $\lambda_9(S) = \lambda_0(\dots S_{-27-2n}(21)^n\,1111\,\omega_1\,\omega_1\,\omega_1^*\,\omega_1\,2122\,1212\,(333212)^n\,S_{22+6n}\dots) \\
        \leq \lambda_0(\overline{21}\,1111\,\omega_1\,\omega_1\,\omega_1^*\,\omega_1\,2122\,1212\,(333212)^n\,3\dots)\leq \lambda_0(\overline{21}\,1111\,\omega_1\,\omega_1\,\omega_1^*\,\omega_1\,2122\,1212\,\overline{333212})+\varepsilon(n)$.

        By direct computation as before we obtain:
        $$\lambda_9(S)\leq\lambda_0(\overline{21}\,1111\,\omega_1\,\omega_1\,\omega_1^*\,\omega_1\,2122\,1212\,\overline{333212})+\varepsilon(n) < J.$$
        \item $\lambda_{18}(S) = \lambda_0(\dots2111\,\omega_1^*\,2122\dots) <j_0-10^{-6}$ according to the point $6$ of Lemma \ref{locuniq2}.
        \item Lastly, $\fa 0\leq k < n $: $$\max_{i\in \{-1,0,1\}}\lambda_{6k+32+ i}(S) = \lambda_{6k+31}(S) = \lambda_0(\dots2123^*33212\dots)<j_0 -10^{-3}.$$
    \end{enumerate}
So $\lambda_0(S) \in M$ and:
$$X = \left\{\lambda_0\left(S_l^T\,(21)^n\,1111\,\omega_1\,\omega_1^*\,\omega_1\,\omega_1\,2122\,1212(333212)^n\,S_r\right) \: / \:(S_l,S_r) \in \{1,2\}^{\N}\right\}\subset M \cap [j_1,j_1+\alpha).$$

We define two dynamical Cantor sets:

$$A_n = \left\{ [3,2,1,1,1,(\omega_1)^2, 2,1,2,2,(1,2)^2, (3,3,3,2,1,2)^n, S_r]\: / \:S_r\in \{1,2\}^{\N} \right\},$$
$$B_n = \left\{ [0,3,2,1,2,\omega_1^T,1,1,1,1,(1,2)^n, S_l]\: /\: S_l\in \{1,2\}^{\N} \right\}.$$
So $X = A_n + B_n$, where $A_n, B_n$ are two sets diffeomorphic to:
$$C(2) := \left\{ [0,a_1,a_2,\dots] \:/ \:\fa n \in \N \:a_n\in \{1,2\} \right\} \: \text{with} \: \: HD(C(2)) >0.5,$$
where $C(2)$ is a regular Cantor set of class $\mathcal{C}^2$ non-essentially affine (see \cite[Proposition 1]{geometricproperties}. Thus we have:
$$HD(A_n) = HD(B_n)> 0.5$$

So according to the Moreira’s dimension formula \cite{CartesianProduct}:

$$HD(X) = \min \left\{1, HD(A_n) + HD(B_n\,)\:\right\} =1,$$
and therefore:
$$1 \geq HD(M\cap(j_1,j_1+\alpha)) \geq HD(X)  = 1.$$
\end{proof}

\begin{corollary}\label{dimloc}
    We have:
    $$\fa \alpha > 0, \:HD(L\cap (j_1,j_1+\alpha) ) =1.$$
\end{corollary}
\begin{proof}
    It is a consequence from:
    $$HD(M\setminus L) < 1.$$
    Indeed:
    $$M\cap (j_1,j_1+\alpha)= L\cap (j_1,j_1+\alpha) \sqcup (M\setminus L)\cap (j_1,j_1+\alpha),$$
    and:
    $$1 = HD(M\cap (j_1,j_1+\alpha) ) = max\left\{HD(L\cap (j_1,j_1+\alpha)), HD((M\setminus L)\cap (j_1,j_1+\alpha))\right\}.$$
    So necessarily, we must have:
    $$HD (L\cap (j_1,j_1+\alpha)) = 1.$$
\end{proof}

\section{Portion of $M\setminus L$ in the vicinity of $\lambda_0(\overline{12333^*2112})$}
For this section, we denote $\omega_2 = 123332112$ and $\omega_2^* = 12333^*2112$, with a asterisk when we refer to it with the position zero.

\subsection{Local uniqueness}

We call $\Z^- = \rrbracket -\infty, 0 \rrbracket$.
\begin{lemma}\label{fw1}(Forbidden words):\\
Let $S \in \{1,2,3\}^{\Z}$ be a bi-infinite sequence.
\begin{enumerate}
    \item If $S = 13^*$, then $\lambda_0(S)> j_0+10^{-1}$.
    \item If $S = 23^*2$, then $\lambda_0(S)> j_0+10^{-2}$.
    \item If $S = 33^*23$, then $\lambda_0(S)> j_0+10^{-2}$.
    \item If $S = 233^*22 \:\text{or} \:333^*22 $, then $\lambda_0(S)> j_0+10^{-3}$.
    \item If $S = 1\,\omega_2^*\,12$, then $\lambda_0(S)> j_0+10^{-5}$.
    \item If $S = 32\,\omega_2^*\,12,$ then $\lambda_0(S)> j_0+10^{-5}$.
    \item If $S = 12\,\omega_2^*\,121, 12\,\omega_2^*\,122, 22\,\omega_2^*\,121, 22\,\omega_2^*\,122$, then $\lambda_0(S)> j_0+10^{-6}$.
    \item If $S = 122\,\omega_2^*\,12,
    222\,\omega_2^*\,12$, then $\lambda_0(S)> j_0+10^{-6}$.
    \item If $S = 3322\,\omega_2^*\,123, 322\,\omega_2^*\,1233$, then $\lambda_0(S)> j_0+10^{-7}$.
    \item If $S = 12112\,\omega_2^*\,1233, 22112\,\omega_2^*\,1233, 2112\,\omega_2^*\,12332$, then $\lambda_0(S)> j_0+10^{-8}$.
    \item If $S = 332112\,\omega_2^*\,123333, 2332112\,\omega_2^*\,1233321, 3332112\,\omega_2^*\,12333212, \\
    23332112\,\omega_2^*\,123332111, 33332112\,\omega_2^*\,123332111$ then $\lambda_0(S)> j_0+10^{-10}$.
    \item If $S = 33332112\,\omega_2^*\,\omega_2$, then $\lambda_0(S) >j_0+10^{-10}$.
     \item If $S = \omega_2\,\omega_2^*\,\omega_2\,2, \omega_2\,\omega_2^*\,\omega_2\,3 $, then $\lambda_0(S)> j_0+9\times10^{-10}$.
     \item If $S = 2\,\omega_2\,\omega_2^*\,\omega_2\,11$, then $\lambda_0(S) >j_0+4\times10^{-11}$.
     \item If $S = 212\,\omega_2\,\omega_2^*\,\omega_2\,123$, then $\lambda_0(S) >j_0+1.82 \times 10^{-12}$.
     \item If $S = 1112\,\omega_2\,\omega_2^*\,\omega_2\,123$, then $\lambda_0(S) >j_0+1.09\times10^{-13}$.
\end{enumerate}
\end{lemma}

\begin{corollary}\label{cor2}
Let $S\in \{1,2,3\}^{\Z}$ such that $m(S) \leq j_0+10^{-4}$. Then, the words $322,223$ and $323$ are forbidden.
\end{corollary}
\begin{proof}
    If $S$ contains $323$, then since according to the Lemma \ref{fw1}, the words $13, 232$ and $3323$ are forbidden, we can't extend the word $323$ to the left without making a forbidden word appear.
    
    If $S = \dots322\dots$, then according to the Lemma \ref{fw1}, since $13$ and $232$ are forbidden, the word $322$ must extend to the left in such way: $S = \dots3322\dots$. However, since the words $13,23322$ and $33322$ are forbidden, we can't extend the word $3322$ without having a contradiction.
    The same reasoning gives that $223$ is forbidden.    
    
\end{proof}

\begin{lemma}\label{locuniq_1}(Local uniqueness):\\
Let $S \in \{1,2,3\}^{\Z}$ be such that $\lambda_0(S) \leq j_0 +10^{-7}$ then $S$ or $S^T$ must take one of these forms;
\begin{enumerate}
    \item $S = 1^*,2^*$ or $33^*3$ and $\lambda_0(S) <  j_0-10^{-2}$,
    \item $S = 233^*21$ or $333^*212$ and $\lambda_0(S) <  j_0-10^{-3}$,
    \item $S=2333^*2111, 3333^*2111$ or $3333^*2112$ and $\lambda_0(S) <  j_0-10^{-4}$,
    \item $S=\omega_2^*\,2$ or $\omega_2^*\,3$ and $\lambda_0(S) <  j_0-10^{-3}$,
    \item $S=1\,\omega_2^*\,11, 2\,\omega_2^*\,11, 212\,\omega_2^*\,123$  and $\lambda_0(S) <  j_0-10^{-6}$,
    \item $S=1112\,\omega_2^*\,123$  and $\lambda_0(S) <  j_0-10^{-7}$,
    \item $S=2112\,\omega_2^*\,1233$.
\end{enumerate}

\end{lemma}
As a consequence, if $S$ verifies $|\lambda_0(S)-j_0|\leq 10^{-7}$ then $S$ must be of the form: $$S=2112\,\omega_2^*\,1233.$$

\begin{proof}
Let $S \in \{1,2,3\}^{\Z}$ such that $m(S) \leq j_0+10^{-7}$. Then, according to Lemma \ref{fw1}, all words until item $9$ are forbidden.
If $S = 1^*$ or $2^*$, then using the bounding rules, $\lambda_0(S) <j_0-10^{-2}$.

If $S = 3^*$, then, using the  forbidden words $13, 232, 31$, we must have:
    $$S = 23^*3, 33^*2 \:\:\text{or}\:\: \left\{ S= 33^*3\:\:\text{and}\:\: \lambda_0(S) < j_0-10^{-2}\right\}.$$
If $S = 33^*2$, then, using the  forbidden words $13, 323$ and $322$, we must have:
    $$S = 333^*21 \:\:\text{or}\:\:\left\{ S= 233^*21\:\:\text{and}\:\: \lambda_0(S) < j_0-10^{-3}\right\}.$$
If $S = 333^*21$, then, using the  forbidden words $13$, we must have:
    $$S = 2333^*211, 3333^*211 \:\:\text{or}\:\:\left\{ S= 333^*212\:\:\text{and}\:\: \lambda_0(S) < j_0-10^{-3}\right\}.$$
If $S = 3333^*211$, then, using the  forbidden words $13$, we must have:
    $$S = 3333^*2111, 3333^*2112 \:\:\text{and}\:\: \lambda_0(S) < j_0-10^{-3}.$$
If $S = 2333^*211$, then, using the  forbidden words $13, 323, 223$, we must have:
    $$S = \omega_2^*\:\:\text{or}\:\:\left\{ S= 2333^*2111\:\:\text{and}\:\: \lambda_0(S) < j_0-10^{-4}\right\}.$$
If $S = \omega_2^*$, then, using the  forbidden words $31$, we must have:
   $$S = 1\,\omega_2^*\,1, 2\,\omega_2^*\,1\:\:\text{or}\:\:\left\{ S= \omega_2^*\,2, \omega_2^*\,3\:\:\text{and}\:\: \lambda_0(S) < j_0-10^{-3}\right\}.$$
If $S = 1\,\omega_2^*\,1$, then, using the  forbidden words $31, 13, 1\,\omega_2\,12$, we must have:
   $$S = 1\,\omega_2^*\,11, \:\:\text{and}\:\: \lambda_0(S) < j_0-10^{-6}.$$
If $S = 2\,\omega_2^*\,1$, then, using the  forbidden words $13, 32\,\omega_2\,12$, we must have:
   $$S = 12\,\omega_2^*\,12, 22\,\omega_2^*\,12\:\:\text{or}\:\:\left\{ S= 2\,\omega_2^*\,11\:\:\text{and}\:\: \lambda_0(S) < j_0-10^{-6}\right\}.$$
If $S = 22\,\omega_2^*\,12$, then, using the  forbidden words $22\,\omega_2\,121, 122\,\omega_2\,12, 222\,\omega_2\,12, 22\,\omega_2\,122$, we must have:
   $$S = 322\,\omega_2^*\,123,$$
which can no longer be extended, since the words $13, 31, 232$ and $3322\,\omega_2\,123$ are forbidden. So the word $22\,\omega_2\,12$ is forbidden. If $S = 12\,\omega_2^*\,12$, then, using the  forbidden words $13, 12\,\omega_2\,122, 12\,\omega_2\,121$, we must have:
   $$S = 112\,\omega_2^*\,123 \:\:\text{or}\:\:\left\{ S= 212\,\omega_2^*\,123\:\:\text{and}\:\: \lambda_0(S) < j_0-10^{-6}\right\}.$$
If $S = 112\,\omega_2^*\,123$, then, using the  forbidden words $13, 232$, we must have:
   $$S = 2112\,\omega_2^*\,1233 \:\:\text{or}\:\:\left\{ S= 1112\,\omega_2^*\,123\:\:\text{and}\:\: \lambda_0(S) < j_0-10^{-7}\right\},$$
which proves the property of local-uniqueness.
\end{proof}

Again, we can summarize the development of the sequence around the vicinity of $j_0$ with a tree of possibilities.
\subsection{Self-replication}
 
\subsubsection{Simplification of forbidden words}
\begin{lemma}\label{fw2}
    Let $S \in \{1,2,3\}^{\Z}$ be a bi-infinite sequence.
    \begin{enumerate}
        \item If $S = 33^*2$, then $m(S) > j_0 +10^{-2}$.
        \item If $S = 3^*22$, then $m(S) > j_0 +10^{-3}$.
        \item If $S = 22\,\omega_2^*\,12$, then $m(S) > j_0 +10^{-7}$.
        \item  If $S = 33332112\,\omega_2^*\,12333211$, then $m(S) > j_0 +9\times 10^{-11}$.
    \end{enumerate}
\begin{proof}
    If $S = 3^*23$ then using the forbidden words $13, 232$ and $3323$, we see, according to point 1 and 2 of Lemma \ref{fw1} that we must have:
    $$m(S)>j_0+10^{-2}.$$
    If $S = 3^*22$, then using the forbidden words $13$ and $232$ (otherwise, we also get $m(S)>j_0+10^{-2})$, we must have $S = 33^*22$, then according to point 3 of Lemma \ref{fw1}, we must have:
    $$m(S) >j_0+10^{-2}.$$
    If $S = 22\,\omega_2^*\,12$, then according to point 8 of the lemma \ref{fw1}, if $S_{-7} \in \{1,2\}$ or $S_{7} \in \{1,2\}$, then $\lambda_0(S)>j_0 +10^{-6}$.
    Otherwise, $S = 322\,\omega_2^*\,123$ and using the forbidden words $3322\,\omega_2\,123, 232, 13, 31$ (otherwise, we also get $m(S)>j_0+10^{-7})$, we must have: 
    $$\inf_{S_g,S_d \in \{1,2,3\}}\sup_{n\in \Z} \lambda_n(\dots S_g\,322\,\omega_2^*\,123\,S_d\dots)>j_0+10^{-7}.$$
    If $S = 33332112\,\omega_2^*\,12333211$, then using the forbidden words $33332112\,\omega_2\,12332111, 33332112\,\omega_2\,, 13$, according to Lemma \ref{fw1}, we must have:
    $$\inf_{S_g,S_d \in \{1,2,3\}}\:\sup_{n\in \Z} \lambda_n(\dots S_g\,33332112\,\omega_2^*\,12333211\,S_d\dots)>j_0+10^{-11}.$$
    If $S = 33332112\,\omega_2^*\,12333211$, then we have:
    $$\inf_{S_g,S_d \in \{1,2,3\}^3}\sup_{n\in \Z}\lambda_n(\dots S_g\,33332112\,\omega_2^*\,12333211\,S_d\dots) \geq \lambda_0( \overline{31}\, 123\,33332112\,\omega_2^*\,12333211\,212\,\overline{13}),$$
    with:
   $$ \lambda_0( \overline{31}\, 123\,33332112\,\omega_2^*\,12333211\,212\,\overline{13}) \geq j_0 + 9\times 10^{-11}.$$
\end{proof}   
    
\end{lemma}

Using the "new" forbidden words, we now call $F_{tot}$ the set of every forbidden word (and its transposes), in the developping sequence of $\overline{123332112}$. We have:
$$F= \{13,232 ,
 322,323,
  1 \,\omega_2\,12 ,
  22 \,\omega_2\,12 ,
  32 \,\omega_2\,12 ,
  12 \,\omega_2\,121 ,
  12 \,\omega_2\,122 ,$$
  $$22112 \,\omega_2\,1233 ,
  2112 \,\omega_2\,12332 ,
  12112 \,\omega_2\,1233 ,
  332112 \,\omega_2\,123333 ,
  2332112 \,\omega_2\,1233321 ,
  3332112 \,\omega_2\,12333212 ,$$
  $$33332112 \,\omega_2\,12333211 ,
  23332112 \,\omega_2\,123332111 ,
   \,\omega_2\,\omega_2\,\omega_2\,2 ,
   \,\omega_2\,\omega_2\,\omega_2\,3 ,
  2 \,\omega_2\,\omega_2\,\omega_2\,11 ,
  212 \,\omega_2\,\omega_2\,\omega_2\,123 ,
  1112 \,\omega_2\,\omega_2\,\omega_2\,123 \},$$
and:
$$F_{tot} = F\cup F^T.$$
Let us call:
\begin{enumerate}
    \item $w_j = \textcolor{red}{1}112 \,\omega_2\,\omega_2\,\omega_2\,123$,
    \item $w_{j'} =\textcolor{red}{2}12 \,\omega_2\,\omega_2\,\omega_2\,123$,
    \item $w_J =  2 \,\omega_2\,\omega_2\,\omega_2\,1\textcolor{red}{1}$.
\end{enumerate}

\begin{theorem}\label{w_sr_1}(Self-replication)

    Let $S  \in \{1,2,3\}^{\Z}$ such that: $S_{n-8}\dots S_{n+6} = 2112\,\omega_2^*12$ for $n \in \Z$.
    \begin{enumerate}
        \item If the finite word $S_{n-15}\dots S_{n+15}$ doesn't contain any word of $F_{tot}\setminus\{w_j,w_j^T, w_{j'},w_{j'}^T\}$, then $S_{n-15}\dots S_{n+15} = 12\,\omega_2\,\omega_2^*\,\omega_2\, 12$.
        \item If the sequence $S_{n-17}\dots S_{n+15}$ doesn't contain any word of $F_{tot}$, then $S_{n-17}\dots S_{n+15} = 2112\,\omega_2\,\omega_2^*\,\omega_2\, 12$.
    \end{enumerate}
       
\end{theorem}

\begin{proof}
    If $S = 2112\,\omega_2^*\,12$, then, using the forbidden words $12\,\omega_2\,121, 12\,\omega_2\,122$, we must have:
   $$S = 2112\,\omega_2^*\,123. $$
   If $S = 2112\,\omega_2^*\,123$, then, using the forbidden words $31, 232$, we must have:
   $$S = 2112\,\omega_2^*\,1233. $$
   If $S = 2112\,\omega_2^*\,1233$, then, using the forbidden words $31, 22112\,\omega_2\,1233, 2112\,\omega_2\,12332, 12112\,\omega_2\,1233$, we must have:
   $$S = 32112\,\omega_2^*\,12333. $$
   If $S =32112\,\omega_2^*\,12333$, then, using the forbidden words $13, 31, 232, 332112\,\omega_2\,123333$, we must have:
   $$S = 332112\,\omega_2^*\,123332. $$
   If $S =332112\,\omega_2^*\,123332$, then, using the forbidden words $13, 323, 322, 2332112\,\omega_2\,1233321$, we must have:
   $$S = 3332112\,\omega_2^*\,1233321. $$
   If $S =3332112\,\omega_2^*\,1233321$, then, using the forbidden words $13, 3332112\,\omega_2\,12333212, 33332112\,\omega_2\,12333211$ we must have:
   $$S = 23332112\,\omega_2^*\,12333211. $$
   
   If $S =23332112\,\omega_2^*\,12333211 $, then, using the forbidden words $13, 3233, 223333, 23332112\,\omega_2\,123332111$, we must have:
   $$S =\omega_2\,\omega_2^*\,\omega_2.$$
   If $S =\omega_2\,\omega_2^*\,\omega_2 $, then, using the forbidden words $31, 1\,\omega_2\,12, \omega_2\,\omega_2\,\omega_2\,2, \omega_2\,\omega_2\,\omega_2\,3$, we must have:
   $$S =2\,\omega_2\,\omega_2^*\,\omega_2\,1.$$
   If $S =2\,\omega_2\,\omega_2^*\,\omega_2 \,1$, then, using the forbidden words $13, 22\,\omega_2\,12, 32\,\omega_2\,12, 2\,\omega_2\,\omega_2\,\omega_2\,11$, we must have:
   $$S =12\,\omega_2\,\omega_2^*\,\omega_2\,12,$$
   which proves the first point. Again, using the forbidden words, $31, w_{j'} = 212\,\omega_2\,\omega_2\,\omega_2\,123$ we must have:
   $$S =112\,\omega_2\,\omega_2^*\,\omega_2\,12.$$
   If $S =2\,\omega_2\,\omega_2^*\,\omega_2 \,1$, then, using the forbidden words $31, w_j = 1112\,\omega_2\,\omega_2\,\omega_2\,123$, we must have:
   $$S =2112\,\omega_2\,\omega_2^*\,\omega_2\,12,$$
   which concludes the proof.
\end{proof}

Using the Lemma \ref{w_sr_1} repeatedly, we have two following results.
\begin{corollary}\label{w_selfrep}
    Let $S  \in \{1,2,3\}^{\Z}$ such that: $S_{-8}\dots S_{8} = 2112\,\omega_2^*12$.
    \begin{enumerate}
        \item Then, if all the words from $F_{tot}\setminus\{w_j,w_j',w_j^T,w_j'^T\}$ are forbidden in $S$, then we must have:
    $$S = S_l^T12\,\omega_2\,\omega_2^*\,\overline{\omega_2} \:\:\text{with}\:\: S_l\in \{1,2,3\}^{\N}.$$
        \item If all the words from $F_{tot}$ are forbidden in $S$, then we must have:
    $$S = \overline{\omega_2} = \overline{12333^*2112}.$$
    
    \end{enumerate}
\end{corollary}

\begin{lemma}\label{w_j_J}
We have:
\begin{enumerate}
     \item 
     \begin{align*}
     J &= 
     \begin{multlined}[t]
     \min\{\lambda_0(S_l^T \,w_J^* \, S_r)\: /\: (S_r,S_l) \in \{1,2,3\}^{\N}\:\text{and} \:\:
     S_l^T w_JS_r \:\:\text{doesn't contain any words from} \\ F_{tot}\setminus\{w_j,w_{j'},w_J,w_j^T,w_{j'}^T,w_J^T\}\} 
     \end{multlined} \\ 
     &=\lambda_0(\overline{12}\,\omega_2\,\omega_2\,\omega_2^*\,\omega_2\,11\,\overline{12}) \approx j_0 + 5.88429645\times 10^{-11}.
    \end{align*}
    \item  \begin{align*}
    j' &= 
    \begin{multlined}[t]
    \min\big\{ \lambda_0(S_l^T w_{j'}^*S_r) / (S_l,S_r)\in \{1,2,3\}^{\N}\:\:\text{and}\:\:S_l^T w_{j'}S_r \:\:\text{doesn't contain any words from} \\  F_{tot}\setminus\{w_j,w_{j'},w_j^T,w_{j'}^T\}\big\} 
    \end{multlined}\\
    &= \lambda_0(\overline{21}\,12332121\,12333212\,\omega_2\,\omega_2^*\,\overline{\omega_2}) \approx j_0 + 2.2055806\times 10^{-12}.
    \end{align*}

    \item  \begin{align*}
    j &= \min\left\{ \lambda_0(S_l^T w_j^*S_r) / (S_l,S_r)\in \{1,2,3\}^{\N}\:\:\text{and}\:\:S_l^T w_jS_r \:\:\text{doesn't contain any words from}\:\:  F_{tot}\setminus\{w_j,w_j^T\}\right\} \\
    &= \lambda_0(\overline{21}\,1112\,\omega_2\,\omega_2^*\,\overline{\omega_2}) \approx j_0 + 4.77646040\times 10^{-13}.  
    \end{align*}
\end{enumerate}

In particular:
$$j<j'<J.$$

\begin{remark}
    We also have:
    \begin{enumerate}
        \item $j = \lambda_0(\overline{21}\,1112\,\omega_2\,\omega_2^*\,\overline{\omega_2}) = m(\overline{21}\,1112\,\omega_2\,\omega_2^*\,\overline{\omega_2}).$
        \item $j' = \lambda_0(\overline{21}\,12332121\,12333212\,\omega_2\,\omega_2^*\,\overline{\omega_2}) = m(\overline{21}\,12332121\,12333212\,\omega_2\,\omega_2^*\,\overline{\omega_2}).$
        \item  $J = \lambda_0(\overline{12}\,\omega_2\,\omega_2\,\omega_2^*\,\omega_2\,11\,\overline{12}) = m(\overline{12}\,\omega_2\,\omega_2\,\omega_2^*\,\omega_2\,11\,\overline{12}).$
    \end{enumerate}
    So:
    $$(j,j',J)\in M^3.$$

    It is relevant to understand the geometric role of these constants. We will prove later that $j$ is the first element in $M$ coming next after $j_0$ and is the left border of the region of $M \setminus L$ to the right of $j_0$, i.e. $\min(M\setminus L)\cap(j_0,J)=j$. It is also the first element of $M$ in the vicinity of $j_0$, containing a forbidden word, which is $w_j$. The point $j'$ is the first element of $M$ containing another forbidden word, which is $w_{j'}$. Finally, $J$ represents the right border of the region of $M\setminus L$ to the right of $j_0$, and is also the first element of $M$ containing the forbidden word $w_J$.
    
    We will prove later that $J$ is a point in $L'$ since above it, sequences $S\in\{1,2,3\}^\Z$ such that $m(S)=\lambda_0(S)$, can be different from $\overline{\omega_2}$ in both directions (left and right), whereas below $J$ and above $j_0$, it can only be different to $\overline{\omega_2}$ to the left direction. 
\end{remark}

\end{lemma}
\begin{proof}
    \textbf{Proof of $J$:}\\
    
    Let $S = S_l^Tw_JS_r$ such that $S$ doesn't contain any words from  $F_{tot}\setminus\{w_j,w_j^T,w_{j'},w_{j'}^T,w_J,w_J^T\}$.
    Again, we have:
    $$S = \dots S_{-15}\,2\,\omega_2\,\omega_2^*\,\omega_2\,11\,S_{16}\dots$$

    In order to minimise $\lambda_0$, since $13$ is forbidden, the right extension must be:
    $$S_{16}\dots = \overline{12}.$$
    So:
    $$S = \dots S_{-15}\,2\,\omega_2\,\omega_2^*\,\omega_2\,11\,\overline{12}.$$
    Using the forbidden words $22\,\omega_2\,12$ and $32\,\omega_2\,12$, we must have:
    $$S= \dots S_{-16}\,12\,\omega_2\,\omega_2^*\,\omega_2\,11\,\overline{12}.$$
    Since we are minimising $\lambda_0(S)$, we must have:
    $$S = \dots S_{-17}\,112\,\omega_2\,\omega_2^*\,\omega_2\,11\,\overline{12}.$$
    Using the forbidden words $31$, we must have:
    $$S = \dots S_{-18}\,2112\,\omega_2\,\omega_2^*\,\omega_2\,11\,\overline{12}.$$
    Then, we can apply the first point of Theorem \ref{w_selfrep} to $(S_{n-9})_{n\in \Z}$, so we get:
    $$S = \dots S_{-25}\,12\,\omega_2\,\omega_2\,\omega_2^*\,\omega_2\,11\,\overline{12}.$$
    Since we are minimising $\lambda_0(S)$, then we must have:
    $$S = \overline{12}\,\omega_2\,\omega_2\,\omega_2^*\,\omega_2\,11\,\overline{12}.$$
    Therefore:
    $$J = \lambda_0(\overline{12}\,\omega_2\,\omega_2\,\omega_2^*\,\omega_2\,11\,\overline{12}) \approx j_0 + 5.88429645\times 10^{-11}.$$

    \textbf{Proof of $j'$:}\\
    Let $S = S_l^Tw_{j'}S_r$ such that $S$ doesn't contain any words from  $F_{tot}\setminus\{w_j,w_j^T,w_{j'},w_{j'}^T\}$.
    Again, we have:
    $$S = \dots S_{-17}\,212\,\omega_2\,\omega_2^*\,\omega_2\,12\,S_{16}\dots$$
    Then, using Corollary \ref{w_selfrep}, we must have:
    $$S = \dots S_{-17}\,212\,\omega_2\,\omega_2^*\,\overline{\omega_2}.$$
    Then, since we are minimising $\lambda_0(S)$, we must have:
    $$S = \dots S_{-18}\,3212\,\omega_2\,\omega_2^*\,\overline{\omega_2}.$$
    Then, using the forbidden words $13$ and $232$, we must have:
    $$S =\dots S_{-19}\,33212\,\omega_2\,\omega_2^*\,\overline{\omega_2}.$$
    Then, since we are minimising $\lambda_0(S)$, we must have:
    $$S =\dots S_{-20}\,333212\,\omega_2\,\omega_2^*\,\overline{\omega_2}.$$
    Then, using the forbidden words $13$ and by minimising $\lambda_0(S)$, we must have:
    $$S =\dots S_{-21}\,2333212\,\omega_2\,\omega_2^*\,\overline{\omega_2}.$$
    Then, using the forbidden words $323$ and $322$, we must have:
    $$S =\dots S_{-22}\,12333212\,\omega_2\,\omega_2^*\,\overline{\omega_2}.$$
    Then, by minimising $\lambda_0(S)$ we can add $2121$, so we get:
    $$S = \dots S_{-26}\,2121\,12333212\,\omega_2\,\omega_2^*\,\overline{\omega_2} = \dots S_{-26}\,21\,\omega_2^T\,2\,123332112\,\omega_2^*\,\overline{\omega_2}.$$
    Then, since $221\,\omega_2^T\,21$ and $221\,\omega_2^T\,21$ are forbidden, we must have:
    $$S = \dots S_{-27}\,32121\,12333212\,\omega_2\,\omega_2^*\,\overline{\omega_2}.$$
    Then, using the forbidden words $13$ and $232$, we must have:
    $$S = \dots S_{-28}\,332121\,12333212\,\omega_2\,\omega_2^*\,\overline{\omega_2}.$$
    Then, using the forbidden words $13$ and minimising $\lambda_0(S)$, we must have:
    $$S = \dots S_{-29}\,2332121\,12333212\,\omega_2\,\omega_2^*\,\overline{\omega_2}.$$
    Then, using the forbidden words $323$ and $322$, we must have:
    $$S = \dots S_{-30}\,12332121\,12333212\,\omega_2\,\omega_2^*\,\overline{\omega_2}.$$
    Finally, we minimise with the periodic sequence, $\overline{21}$, so we get:
    $$S = \overline{21}\,12332121\,12333212\,\omega_2\,\omega_2^*\,\overline{\omega_2}.$$
    Therefore:
    $$j' = \lambda_0(\overline{21}\,12332121\,12333212\,\omega_2\,\omega_2^*\,\overline{\omega_2}) \approx j_0 + 2.2055806\times 10^{-12}.$$
    \textbf{Proof of $j$:}\\
    Let $S = S_l^Tw_{j}S_r$ such that $S$ doesn't contain any words from  $F_{tot}\setminus\{w_j,w_j^T,w_{j},w_{j}^T\}$.
    Again, we have:
    $$S = \dots S_{-18}\,1112\,\omega_2\,\omega_2^*\,\omega_2\,123\,S_{17}\dots$$
    Then, using Corollary \ref{w_selfrep}, we must have:
    $$S = \dots S_{-18}\,1112\,\omega_2\,\omega_2^*\,\overline{\omega_2}.$$
    Then, we minimise $\lambda_0(S)$ by adding the sequence $\overline{21}$ to the left:
    $$S = \overline{21}\,1112\,\omega_2\,\omega_2^*\,\overline{\omega_2}.$$
    Therefore:
    $$j =  \lambda_0(\overline{21}\,1112\,\omega_2\,\omega_2^*\,\overline{\omega_2}) \approx j_0 + 4.77646040\times 10^{-13}.$$  
\end{proof}

\begin{prop}\label{w_remark fw}
    Let $S\in\{1,2,3\}^\Z$ be a sequence such that $\lambda_0(S) = m(S) = \sup_{n \in \Z} \lambda_n(S)$. According to Lemma \ref{fw1} and \ref{w_j_J} we have:

\begin{enumerate}
    \item If $\lambda_0(S) \in M\cap (j_0,J+\alpha)$, for $\alpha >0$ small enough, then, at least one of the subwords $w_j,w_j^T,w_{j'}, w_{j'}^T, w_J,w_J^T$ must appear in $S$.
    \item  If $\lambda_0(S) \in M\cap(j_0,J)$, no forbidden words from $F_{tot}\setminus \{w_j,w_{j'},w_j^T,w_{j'}^T\}$ are allowed to appear in $S$.
    \item  If $\lambda_0(S) \in M\cap(j_0,j')$, no forbidden words from $F_{tot}\setminus \{w_j,w_j^T\}$ are allowed to appear in $S$.
    \item If $\lambda_0(S) \in M\cap(j_0,j)$, no forbidden words from $F_{tot}$ are allowed to appear in $S$.
    
\end{enumerate}

\end{prop}

\begin{proof}

    Let $m\in M$ such that $m>j_0$ and $S\in \{1,2,3\}^{\Z}$ a sequence such that: $$m=\lambda_0(S) = \sup_{n\in \Z}\lambda_n(S).$$
    If $J+\alpha>m>j_0$, with $\alpha< 9\times10^{-10}$, then, according to the Lemma \ref{fw1} and Corollary \ref{w_selfrep}, one of the words $w_j,w_j^{T}, w_{j'}, w_{j'}^T, w_J$ and $w_J^{T}$ must be in the sequence $S$.

    Let assume that $m \in (j_0,J)$. By contradiction, suppose that there exists $w_f\in F_{tot}\setminus \{w_j,w_j^T,w_{j'}, w_{j'}^T\}$ such that $S =  S_l^Tw_fS_r$. If $w_f \in F_{tot}\setminus\{w_j,w_j^T,w_{j'}, w_{j'}^T,w_J,w_J^T\}$, then according to Lemma \ref{fw1} and \ref{fw2}, we must have $\lambda_N(S) >j_0+9 \times 10^{-11}$ for some $N \in \Z$, which is impossible. \textbf{So $S$ doesn't contain any word from $F_{tot}\setminus\{w_j,w_j^T,w_{j'}, w_{j'}^T,w_J,w_J^T\}$}. If $S$ contains $w_J$, then we have $S = S_l^Tw_JS_r$. By definition of $J$, we must have $\lambda_N(S) \geq J$ for some $N \in \Z$ and therefore, $m\geq J$, which is impossible. The same is true with $w_J^T$. \textbf{So $S$ doesn't contain any subwords from $F_{tot}\setminus\{w_j,w_j^T, w_{j'}, w_{j'}^T\}$.}
    
    Let assume that $m \in (j_0,j')$.  Then according to above, $S$ doesn't contain any words from $F_{tot}\setminus\{w_j,w_j^T, w_{j'}, w_{j'}^T\}$.
    If $S = S_l^Tw_{j'}S_r$ then by definition of $j'$, we must have $\lambda_N(S)\geq j'$ and therefore: $m(S) \geq j'$. Again, the same is true for $w_{j'}^T$. Hence, the sequence $S$ doesn't contain any words from $F_{tot}\setminus\{w_j,w_j^T\}$.
    
    Finally, let assume that $m \in (j_0,j)$. Then according to above, $S$ doesn't contain any words from $F_{tot}\setminus\{w_j,w_j^T\}$.
    If $S = S_l^Tw_{j}S_r$, then by definition of $j$, we must have $\lambda_N(S)\geq j$ and therefore: $m(S) \geq j$. Again, the same is true for $w_{j}^T$. Hence, the sequence $S$ doesn't contain any words from $F_{tot}$.
\end{proof}

The main consequence that follows from all of this:
\begin{theorem}\label{principal}
    
    We have:
    \begin{enumerate}
        \item $M\cap(j_0,j) = \emptyset$.
        \item $L\cap(j_0,J) = \emptyset$.
    \end{enumerate}
    
\end{theorem}

\begin{proof}
    Let assume by contradiction that there exist $m\in M\cap(j_0,j)$. Then we can find a sequence $S$ such that $m = \sup_{n\in \Z}\lambda_n(S) = m(S) = \lambda_0(S)$.\:
    Then, we have $j_0<\lambda_0(S) < j< j_0 +   10^{-12}$.
    
    So according to the Lemma \ref{locuniq_1}, we must have:
    $S_{-8}\dots S_{8} = 2112\,\omega_2^*\,1233$. In addition, according to the Proposition \ref{w_remark fw}  the words from the set $F_{tot}$ are forbidden in $S$. 
    Then, because of the Corollary \ref{w_selfrep}, the sequence $S$ must extend in such way:
    $$S = \overline{\omega_2},$$
    and so $m(S) = j_0$, a contradiction.
    
    Let assume by contradiction that there exists $l\in L \cap (j_0,J)$. We use the fact that the Markov values of periodic sequences is dense in $L$ (see \cite[Theorem 2, Chapter 3]{Cusick-Flahive}). Therefore, $\exists (l_n)_{n\in \N} \in L^{\N}$ such that: 
    $$\fa n \in \N, \exists \,\sigma^{(n)} \in \N^{(\N)}, l_n = m\left(\overline{\sigma^{(n)}}\right)\:\:\text{and}\:\:
    \lim_{n\rightarrow \infty} l_n = l.$$
    Let us write $\fa n \in \N, S^{(n)} = \overline{\sigma^{(n)}}$. We can assume (even if it means taking $n \geq n_0$ with $n_0$ big enough), $\fa n \in \N, S^{(n)} \in \{1,2,3\}^{\Z}$ and that $l_n<J$.
    So we have $\lambda_0(S^{(n)}) = l_n <J$. Using  Lemma \ref{locuniq_1}, we have that $S_{-8}^{(n)}\dots S_8^{(n)} = 2112\,\omega_2^*\,1233$ (or its transpose).
    Since $l_n = \lambda_0(S^{(n)}) \in M\cap(j_0,J)$, according to the Proposition \ref{w_remark fw}, no forbidden words from $F_{tot}\setminus \{w_j,w_j^T,w_{j'},w_{j'}^T\}$ are allowed to appear in $S^{(n)}$.
    Therefore, using the Corollary \ref{w_selfrep}, we have:
    $$\fa n \in \N, \:\: S^{(n)} = S_l^{(n)}12\,\omega_2\,\omega_2^*\,\overline{\omega_2}\:\:\text{with}\:\:S_l^{(n)} \in \{1,2,3\}^{\Z^-}.$$
    Hence $\fa n \in \N,\:\: \overline{\sigma^{(n)}} = S_l^{(n)}12\,\omega_2\,\omega_2^*\,\overline{\omega_2}$ and necessarily $ \overline{\sigma^{(n)}} = S^{(n)} = \overline{\omega_2}$  since $\sigma^{(n)}$ is periodic. So for all $n$ we have $l_n = m(\overline{\omega_2}) = j_0$ and $l = j_0$, which is impossible.

\end{proof}

\subsection{A portion of $M\setminus L$}

Now, we can characterize the set $(M\setminus L)\cap(j_0,J)$.

\begin{corollary}\label{w_carac1}
    Let $m \in M\cap(j,J)$ and $S \in \{1,2,3\}^{\Z}$ be a sequence such that $m = \lambda_0(S)$. Then we have:
    $$S = S_l\,t\,12\,\omega_2\,\omega_2^*\,\overline{\omega_2}, $$
    
    with $S_l\in \{1,2,3\}^{\Z^-}$ and $t\in \{2,11\}$ such that:
    
    \begin{enumerate}
        \item The sequence $S$ cannot contain any words from $F_{tot}\setminus\{w_j,w_j^T,w_{j'},w_{j'}^T\}$.
        \item $S_l\,t\,12$ cannot contain the word $2112\,\omega_2\,12$.
        \item If $\exists n \in\Z^-$ such that $S_{n-6}^T\dots S_{n+8}^T = 21\,\omega_2^T\,2112$, then we have $S = \overline{\omega_2^T}\,21\,S_{n+16}\dots$.
    \end{enumerate}
   
 \end{corollary}
\begin{proof}
    Let $m(S) = \lambda_0(S) \in (j_0,J)$. According to the Lemma \ref{locuniq_1}, we have:
    $$S_{-8}\dots S_6 = 2112\,\omega_2^*\,12.$$
    Since $m(S) < J$, the sequence $S$ doesn't contain any words from $F_{tot}\setminus\{w_j,w_{j'},w_j^T,w_{j'}^T\}$. Then, because of the Theorem \ref{w_selfrep}, we have:
    $$S = \dots S_{-16}12\,\omega_2\,\omega_2^*\,\overline{\omega_2}.$$
    Since the word $31$ is forbidden, we have $S_{-16}\in \{1,2\}$. 
    
    If $S_{-16}=1$, then using the forbidden word $31$, we have $S_{-17}\in \{1,2\}$. If $S_{-17} = 2$, then because of Theorem \ref{w_selfrep}, applied to $(S_{n-9})_{n\in \Z}$ we must have:
        $$S = \dots12\,\omega_2\,\omega_2\,\omega_2^*\overline{\omega_2},$$
        and $\lambda_0(S)<j_0+10^{-19}<j$, which is impossible. So $S_{-17} = 1$ (and so $t = 11$).
        
    Therefore $S = S_l\,t\,12\,\omega_2\,\omega_2^*\,\overline{\omega_2}$ with $S_l\in \{1,2,3\}^{\Z^-}$ and $t\in \{2,11\}$. If $S_l\,t\,12$ contains the word $2112\,\omega_2\,12$, then because of self-replication to the right, we have no uniqueness in the writing of $S$. So $S_l\,t12$ does not contain the word $2112\,\omega_2\,12$.
        If $\exists n \in \Z^-$ such that:
        $$S_{n-6}\dots S_{n+8} = 21\,\omega_2^T\,2112.$$
        Then, since the words from $F_{tot}\setminus\{w_j,w_j^T,w_{j'},w_{j'}^T\}$ are forbidden, using self replication applied to the sequence $S$ to the left, we have:
        $$S = \overline{\omega_2^T}\,21\,S_{n+16}\dots$$
\end{proof}

Before making a global description of the set $(M\setminus L) \cap (j_0,J)$, we can simplify a little bit the set of forbidden words.
\begin{lemma}\label{simplif}
    Let us call: 
    $$F_0 = \left\{13,232,322,323,1 \,\omega_2\,12 ,
  22 \,\omega_2\,12 ,
  32 \,\omega_2\,12 ,
  12 \,\omega_2\,121 ,
  12 \,\omega_2\,122 , 2112\,\omega_2\,12\right\}.$$
  Then, if $S\in \{1,2,3\}^{\Z}$ is a sequence such that the set $F_0$ is forbidden in $S$, then the words $\omega_2\,121$ and $\omega_2\,122$ are forbidden as well.
\end{lemma}
\begin{proof}
    Indeed, if $S = \dots \omega_2\,121\dots$, then using the forbidden words $31$ and $1\,\omega_2\,12$, then we mjust have:
    $$S = \dots2\,\omega_2\,121\dots$$
    But we can no longer extend $S$ to the left since the words $32\,\omega_2\,12,22\,\omega_2\,12$ and $12\,\omega_2\,121$ are forbidden. The same is true for the word $w\,122$.
\end{proof}
We can now make a global description of the set $(M\setminus L) \cap (j_0,J)$.
\begin{theorem}\label{description globale M/L second region}
    Let us call $F$ the following set:
    $$F = \left\{13,232,322,323,1 \,\omega_2\,12 ,
  22 \,\omega_2\,12 ,
  32 \,\omega_2\,12 ,
  \,\omega_2\,121 ,
  \,\omega_2\,122 , 2112\,\omega_2\,12\right\},$$
  and:
  $$\Tilde{F} = F\cup F^T.$$
  
  We have:
  $$(M\setminus L)\cap(j_0,J) = C \cup D \cup X,$$
  with:
  $$C = \left\{\lambda_0(S\,t\,12\,\omega_2\,\omega_2^*\,\overline{\omega_2}) / \:\:(S,t) \in \{1,2,3\}^{\N}\times \{2,11\}  \text{and}\:\: S\,t\,12\:\:\text{doesn't contain any words from} \:\:\Tilde{F}\right\},$$
  $$D = \{\lambda_0(\overline{\omega_2^T}\,21\,s\,12\,\omega_2\,\omega_2^*\,\overline{\omega_2})\:\:/s\in \{1,2,3\}^{N}, N\geq 1\:\text{such that}\:\:(s_1s_2,s_{N-1}s_N) \in \{11,2\}^2 , [0,s^T]\leq [0,s]$$
  $$
  \text{ and }\:\: 1\,\omega_2^T\,\omega_2^T\,21\,s\,12\,\omega_2\,\omega_2\,1\:\:\text{ doesn't contain any words from }\Tilde{F}\},$$
  and:
  $$X = \left\{\lambda_0(\overline{\omega_2^T}\,2\,\omega_2\,\omega_2^*\,\overline{\omega_2}), \lambda_0(\overline{\omega_2^T}\,212\,\omega_2\,\omega_2^*\,\overline{\omega_2}), \lambda_0(\overline{\omega_2^T}\,21112\,\omega_2\,\omega_2^*\,\overline{\omega_2})\right\}.$$
  
\end{theorem}
\begin{proof}
    Let $m\in (M\setminus L)\cap(j_0,J)$ and $S\in\{1,2,3,\}^\Z$ a sequence such that $\lambda_0(S) = m$. Then according to point $1$ and $2$ of Corollary \ref{w_carac1}, we have:
    $$S = S_l\,t\,12\,\omega_2\,\omega_2^*\,\overline{\omega_2},$$
    with $t \in \{2,11\}$ and $S$ such that the sequence $S$ cannot contain any words from $F_{tot}\setminus\{w_j,w_j^T,w_{j'},w_{j'}^T\}\cup \{21\,\omega_2\,2112\}$. If in addition, the sequence $S$ doesn't contain the word $21\,\omega_2^T\,2112$, then we can simplify the set of forbidden words, since every word from $F_{tot}$ containing the subword $21\,\omega_2\,2112$ can no longer appear. Using the previous Lemma \ref{simplif}, the new set of forbidden words of $S_l\,t\,12$ is exactly $\Tilde{F}$ and therefore, $m = \lambda_0(S) \in C$.

If the sequence $S$ contains the word $21\,\omega_2^T\,2112$, then, let us define $N = \max\{ n\in\Z \:/\:S_{n-6}\dots S_{n+8} = 21\,\omega_2^T\,2112\}$. By definition of $N$:
    $$S = \dots21\,\omega_2^T\,2112\,S_{N+14}\dots$$
    Then according to point $3$ of Corollary \ref{w_carac1}, we have:
    $$S = \overline{\omega_2^T}\,21\,S_{N+21}\dots$$
    And by definition of $N$, we must have $S_{N+21}S_{N+22}\in \{11,2\}$. Now we have two writings of the sequence $S$. So $S$ can be of the following forms:
    \begin{enumerate}
        \item $S = \overline{\omega_2^T}\,2\,\omega_2\,\omega_2^*\,\overline{\omega_2}, \: (N=-33)$,
        \item $S = \overline{\omega_2^T}\,212\,\omega_2\,\omega_2^*\,\overline{\omega_2}, \: (N=-35)$,
        \item $S = \overline{\omega_2^T}\,21112\,\omega_2\,\omega_2^*\,\overline{\omega_2}, \: (N=-37)$,
        \item $S = \overline{\omega_2^T}\,21\,s\,12\,\omega_2\,\omega_2^*\,\overline{\omega_2}, \: (N\leq-37)$.
    \end{enumerate}
    The first three cases correspond to the set $X$. In the last case, $s$ must verify:
    $$\lambda_0(\overline{\omega_2^T}\,21\,s\,12\,\omega_2\,\omega_2^*\,\overline{\omega_2})\geq \lambda_{-31-N}(\overline{\omega_2^T}\,21\,s\,12\,\omega_2\,\omega_2^*\,\overline{\omega_2})= \lambda_0(\overline{\omega_2^T}\,(\omega_2^{T})^*\,\omega_2^T\,21\,s\,12\,\overline{\omega_2}) = \lambda_0(\overline{\omega_2^T}\,21\,s^T\,12\,\omega_2\,\omega_2^*\,\overline{\omega_2}), $$
    which is equivalent to:
    $$[0, s^T] \leq [0,s].$$
    By hypothesis on $N$, the finite sequence  $1\,\omega_2^T\,\omega_2^T\,s\,1\,2\,\omega_2\,\omega_2\,1$ doesn't contain the word $21\,\omega_2^T\,2112$, so like above, we must have that the finite sequence $1\,\omega_2^T\,\omega_2^T\,s\,1\,2\,\omega_2\,\omega_2\,1$ doesn't contain any words from $\Tilde{F}$. In conclusion, we have:
    $$(M\setminus L) \cap(j_0,J) \subset C\cup D \cup X.$$
    Now, we show the reverse inclusion. Let $S_l \in \{1,2,3\}^{\Z^-}$ and $t\in \{2,11\}$ such that $S_l\,t\,12$ doesn't contain any words from $\Tilde{F}$. Let us show that $\lambda_0(S_l\,t\,12\,\omega_2\,\omega_2^*\,\overline{\omega_2}) = m(S_l\,t\,12\,\omega_2\,\omega_2^*\,\overline{\omega_2})$.
    Since we also have $\lambda_0(S_l\,t\,12\,\omega_2\,\omega_2^*\,\overline{\omega_2}) \in (j_0,J)$, it will follow that $\lambda_0(S_l\,t\,12\,\omega_2\,\omega_2^*\,\overline{\omega_2}) \in (M\setminus L) \cap(j_0,J)$.
    We call $S = S_l\,t\,12\,\omega_2\,\omega_2^*\,\overline{\omega_2}$.
    We have:
    \begin{enumerate}
        \item $\fa k \in \N^*, \lambda_{9k}(S)\leq j_0+10^{-19}<j$.
        \item $\lambda_0(S) \geq j$.
        \item $\lambda_{-9}(S) <j_0-10^{-7}$ according to point $5$ and $6$ of Lemma \ref{locuniq_1}.
        \item $\fa k \leq -9$, if we develop the sequence around the index $k$ according to the rules that all the words from $\Tilde{F}$ are forbidden, then we will get that $S_{k-i_g}\dots S_{k+i_d}$ (or its transposes) is a word which belong to the "small" word of Lemma \ref{locuniq_1} and so $\lambda_k(S) <j_0-10^{-7}$.
        So in conclusion, $\lambda_0(S) = m(S)$.
    \end{enumerate}
    We verify by direct computation that $X \subset(M\setminus L)\cap(j_0,J)$. If $m \in D$, then the condition $[0, s^T] \leq [0,s]$ ensures that $m = \lambda_0(\overline{\omega_2^T}\,21\,s\,12\,\omega_2\,\omega_2^*\,\overline{\omega_2}) = m(\overline{\omega_2^T}\,21\,s\,12\,\omega_2\,\omega_2^*\,\overline{\omega_2})$.
    
\end{proof}

\begin{lemma}\label{isolated2}
    The set $D$ previously defined is a set of isolated points of $M\setminus L$.
\end{lemma}
\begin{proof}
    Let $m\in D$ and $s$ be a finite sequence such that $m = \lambda_0(\overline{\omega_2^T}\,21\,s\,12\,\omega_2\,\omega_2^*\,\overline{\omega_2})$. Assume there is a sequence $(m_n)_{n \in\N}\in M^{\N}$ such that:
    $$m = \lim_{n \rightarrow \infty}m_n.$$
    We can assume that we have in addition $\fa n \in \N, m_n\in M\cap (j,J)$.
    Hence, according to the Corollary \ref{w_carac1}, $\fa n \in \N, \exists S_l^{(n)}\in \{1,2,3\}^{\Z^-}$ a sequence such that: 
    $$m_n = \lambda_0(S_l^{(n)}\,t\,12\,\omega_2\,\omega_2^*\,\overline{\omega_2}).$$
    Using the Proposition \ref{continuite lambda},  we can find a integer $n_0$ such that $\fa n\geq n_0$, we have: 
    $$S_l^{(n)} = \Tilde{S_l}^{(n)}\,21\,\omega_2^T\,\omega_2^T\,21\,s.$$
    So the sequence $S_l^{(n)}$ contains the words $21\,\omega_2^T\,2112$ and according to the Theorem \ref{w_carac1}, we must have:
    $$\fa n \geq n_0, \:\: S_l^{(n)} = \overline{\omega_2^T}\,21\,s.$$
    Therefore, $\fa n \geq n_0, \:\: m_n = m$. Hence, $m$ is an isolated point in $(M\setminus L) \cap (j,J)$.
\end{proof}

\begin{lemma}\label{encadrement}
We have:
$$(M\setminus L)\cap(j_0,J) \subset \left[ \lambda_0(\overline{21}\,1112\,\omega_2\,\omega_2^*\,\overline{\omega_2}), \lambda_0(\overline{21}\,212\,\omega_2\,\omega_2^*\,\overline{\omega_2})\right].$$
\end{lemma}
\begin{proof}
    Since we have $M\cap(j_0,j) = \emptyset$ and $L\cap(j_0,J) = \emptyset$, then, we have: 
    $$(M\setminus L)\cap(j_0,J) \subset [j,J).$$
    Then, since we have 
    $\max \, \lambda_0(\dots\,1112\,\omega_2\,\omega_2^*\,\overline{\omega_2}) \leq \min \lambda_0(\dots212\,\omega_2\,\omega_2^*\,\overline{\omega_2})$, then,
    $$\max\,(M\setminus L)\cap(j_0,J) = \max \, \lambda_0(\dots\,212\,\omega_2\,\omega_2^*\,\overline{\omega_2}).$$
    If $S$ is a sequence such that 
    $$\max \,(M\setminus L)\cap(j_0,J) = \lambda_0(S),$$
    then:
    $$S = \dots S_{-17}\,212\,\omega_2\,\omega_2^*\,\overline{\omega_2}.$$
    Since we are maximising the function $\lambda_0$, we must add the sequence $\overline{21}$ to the left, so we get:
    $$S = \overline{21}\,212\,\omega_2\,\omega_2^*\,\overline{\omega_2}.$$
    By computation, we have:
    $$ \lambda_0(\overline{21}\,212\,\omega_2\,\omega_2^*\,\overline{\omega_2})-\lambda_0(\overline{21}\,1112\,\omega_2\,\omega_2^*\,\overline{\omega_2}) \approx 2.409522\times 10^{-12}.$$
\end{proof}

\subsection{The local border of $L$}

\begin{theorem}
    $J \in L' \subset L $ and therefore, $(j_0,J)$ is the largest gap of $L$ near of $j_0$ and we found that $A_2 = (j_0,J)$
\end{theorem}
\begin{proof}
    We have, for all $n \geq 2$
    $$l_n = m\left(\overline{(12)^n\,\omega_2\,\omega_2\,\omega_2^*\,\omega_2\,11\,(12)^n}\right) = l\left(\overline{(12)^n\,\omega_2\,\omega_2\,\omega_2^*\,\omega_2\,11\,(12)^n}\right)=\lambda_0\left(\overline{(12)^n\,\omega_2\,\omega_2\,\omega_2^*\,\omega_2\,11\,(12)^n}\right).$$
    
    So:
    $$\lim_{n \rightarrow \infty} l_n = \lambda_0(\overline{12}\,\omega_2\,\omega_2\,\omega_2^*\,\omega_2\,11\,\overline{12}).$$
    In addition, we have $\fa n \in \N, l_n \in L$. So $J \in L'$ and since the Lagrange spectrum is closed, we finally get that $J\in L$.
\end{proof}

\begin{theorem}
    We have, $\fa \alpha>0, \: HD(M\cap(j_1,j_1+\alpha)) = HD(L\cap(j_1,j_1+\alpha))= 1$.
\end{theorem}

\begin{proof}
    Let $\alpha >0$. Using the continuity of the function $\lambda_0$ according to Proposition \ref{continuite lambda}, there exists $n \in \N$ such that $\fa (S_l,S_r) \in \{1,2\}^{\Z^-}\times \{1,2\}^{\N}$, we have:
    $$\lambda_0(S) = \lambda_0(S_l\,(12)^n\,\omega_2\,\omega_2\,\omega_2^*\,\omega_2\,11\,(12)^n\,S_r) \in [J,J+\alpha).$$
    In addition, we have;
    \begin{enumerate}
        \item $\lambda_9(S)=\lambda_0(\dots2\,\omega_2^*\,11\dots)<j_0-10^{-6}$ according to point 5 from Lemma \ref{locuniq_1}.
        
        \item $\lambda_{-9}(S) \leq \lambda_0(\dots212\,\omega_2\,\omega_2^*\,\omega_2\,\omega_2\dots)$<$j_0+3\times 10^{-12} <J$.
        \item $\lambda_{-18}(S) = \lambda_0(\dots212\,\omega_2^*\,123\dots) <j_0-10^{-6}$ according to point $5$ from Lemma \ref{locuniq_1}.
    \end{enumerate}
    So we get $\lambda_0(S) = m(S)$ and finally, 
    $$X = \left\{\lambda_0(S_l\,(12)^n\,\omega_2\,\omega_2\,\omega_2^*\,\omega_2\,11\,(12)^n\,S_r)\:\:/\:\: (S_l,S_r) \in \{1,2\}^{\Z^-}\times \{1,2\}^{\N} \right\} \subset M \cap [J,J+\alpha).$$
    We define dynamical Cantor sets:
    $$A = \left\{[3,2,1,1,2,w,1,1,(12)^n,S_r]\:\:/\:\: S_r \in \{1,2\}^{\N}\right\},$$
    $$B =  \left\{[3,3,3,2,1,\omega_2^T,\omega_2^T,(21)^n,S_l^T]\:\:/\:\: S_l \in \{1,2\}^{\Z^-}\right\}.$$
    We have $X = A_n + B_n$ and moreover $A_n, B_n$ are diffeomorphic to:
    $$C(2) := \left\{ [0,a_1,a_2,\dots] \:/ \:\fa n \in \N \:a_n\in \{1,2\} \right\} \: \text{with} \: \: HD(C(2)) >0.5,$$
    where $C(2)$ is a regular Cantor set of class $\mathcal{C}^2$, non-essentially affine (see \cite[Proposition 1]{geometricproperties}). We have:
    $$HD(A) = HD(B)> 0.5.$$
    
    So according to the Moreira’s dimension formula \cite{CartesianProduct}:
    
    $$HD(X) = \min \left\{1, HD(A_n) + HD(B_n\,)\:\right\} =1,$$
    and thus:
    $$1 \geq HD(M\cap(j_1,j_1+\alpha)) \geq HD(X)  = 1.$$
\end{proof}

\section{Lower bound on the Hausdorff distance between $M$ and $L$}
In this section, we focus on a portion of the set $M \setminus L$, above the value $b_{\infty} = \lambda_0(\overline{2212^*112}) \approx 3.2930442439$.

We also call the word $w = 2212112$ and $w^* = 2212^*112$.
According to Carlos Gustavo Moreira and Carlos Matheus article \cite{FreimanExample}, we already have many results.
\begin{theorem}
    The largest interval containing $M\setminus L$ near $b_{\infty}$ is $(b_{\infty},B_{\infty})$ with:
    $$B_{\infty} = \lambda_0(\overline{211211212221}\,www^*w\,\overline{112122212112}) \:\:\text{and}\:\: B_{\infty}-b_{\infty}\approx 2.374867\times 10^{-7}.$$
\end{theorem}

They also completely characterized the portion of $M\setminus L$ contained in $(b_{\infty},B_{\infty})$, with the following theorem:
\begin{theorem}\label{central}
    Let $m\in(M\setminus L)\cap (b_{\infty},B_{\infty})$. Then $m = m(B) = \lambda_0(B)$ for a sequence $B\in \{1,2\}^{\Z}$ such that:
    \begin{enumerate}
        \item $B_{-10}\dots B_0B_1\dots B_{7}\dots = ww^*\overline{w}.$
        \item The sequence $(B_n)_{n\leq -11}$ has the following properties:
        \begin{enumerate}
            \item It does not contain the following words: $ 21212 , 121211 ,  112121 , 212111 ,  111212 ,  2w21, $\\
            $ 12w^T2 , 12w22 $.
            \item If $\exists \: n\leq -21$, such that $B_{n-4}\dots B_{n+5} = 2w^T21$, then :
            $\dots B_{n-7}\dots B_{n+10} = \overline{w^T}\,w^T$.
        \end{enumerate}
    \end{enumerate}
\end{theorem}

We now want to estimate the value: 
$$d((M\setminus L)\cap (b_{\infty},B_{\infty}),L) = \sup_{m\in M\cap (b_{\infty},B_{\infty})}d(L,m).$$

Which can be writen like this:
$$d((M\setminus L)\cap (b_{\infty},B_{\infty}),L) = \frac{B_{\infty}-b_{\infty}}{2}-\inf_{m\in M}\left|m - \frac{B_{\infty}+b_{\infty}}{2}\right|.$$
So we only need to focus on the following quantity:
$$\varepsilon_0 = \inf \left\{\left|m-\frac{b_{\infty}+B_{\infty}}{2}\right| \:\: /\:\:m\in M \right\} = \inf \left\{\left|m-\frac{b_{\infty}+B_{\infty}}{2}\right| \:\: /\:\:m\in (M\setminus L) \cap (b_{\infty},B_{\infty}) \right\}.$$
This calculation requires very few steps, thanks to the nice properties of the function $\lambda_0$.

\begin{lemma}\label{encadrement1}
    We have:
    \begin{enumerate}
        \item $\lambda_0(\dots\,2\,ww^*\overline{w}) \leq \lambda_0(\overline{12}\,2\,ww^*\overline{w})\leq \frac{b_{\infty}+B_{\infty}}{2} - 9\times 10^{-8}$.
        \item $\lambda_0(\dots\,1\,ww^*\overline{w}) \geq \lambda_0(\overline{21}\,1\,ww^*\overline{w})\geq \frac{b_{\infty}+B_{\infty}}{2} +2\times 10^{-8}$.
    \end{enumerate}    
\end{lemma}

Because of the Lemma \ref{encadrement1}, to find $\varepsilon_0$, we only need to calculate the following values:
$$m_1= \min\left\{\lambda_0(\,^tB\,1\,ww^*\overline{w}) / \:B\:\:\text{verifies the properties of Theorem \ref{central}}\right\},$$
$$m_2= \max\left\{\lambda_0(\,^tB\,2\,ww^*\overline{w}) / \:B\:\:\text{verifies the properties of Theorem \ref{central}}\right\}.$$
Thus we must have $\varepsilon_0 = \min\{m_1 - \frac{b_{\infty}+B_{\infty}}{2},\frac{b_{\infty}+B_{\infty}}{2}-m_2\}$.
\begin{proof}

Let $S$ be a sequence such that $S = \dots S_{-12}\,1\,ww^*\overline{w}$.

Since we are minimising $\lambda_0(S)$, we must develop $S$ with:
$$S = \dots S_{-16}\,21211\,ww^*\overline{w}.$$
Since the word $121211$ is forbidden, we must have:
$$S = \dots S_{-17}\,221211\,ww^*\overline{w}.$$
Since the word $2w21$ is forbidden, we must have:
$$S = \dots S_{-18}\,1221211\,ww^*\overline{w}.$$

Then we use repeatedly the same argument, since $-18$ and $-12$ are both even and finally, we get:

$$S = \overline{221211}\,ww^*\overline{w},$$
and:
$$m_1  = m(\overline{221211}\,\overline{w}) = \lambda_0(\overline{221211}\,ww^*\,\overline{w}) \approx \frac{b_{\infty}+B_{\infty}}{2} + 2.764903258\times 10^{-8}.$$

\begin{remark}
The point $m_1$ is accumulated on the right by points of $M$.
\end{remark}

The first point of Lemma \ref{encadrement1} ensures that:
$$\varepsilon_0  = m_1 - \frac{b_{\infty}+B_{\infty}}{2}.$$

Therefore:
\begin{equation*}
    \delta_0:= d((M\setminus L)\cap (b_{\infty},B_{\infty}),L)=B_{\infty}-m_1
\end{equation*}
Moreover:
\begin{gather*}
    B_{\infty}-m_1=\lambda_0(\overline{211211212221}\,www^*w\,\overline{112122212112}) - \lambda_0(\overline{221211}\,ww^*\,\overline{w}) \\ 
    =\frac{780369102362985 + 454284412153 \sqrt{151905}}{290742023368320}-\frac{82 \sqrt{87} + 57 \sqrt{18229} + 6931}{4674} \\
    =\frac{272052036746460995 - 3973474319367040 \sqrt{87} - 2762049221999040 \sqrt{18229} + 353887557067187 \sqrt{151905}}{226488036203921280} \\
    \approx 9.1094243388\times 10^{-8}.
\end{gather*}
\end{proof}

\sloppy\printbibliography

\end{document}